\documentclass[11pt]{amsart}  
\usepackage{epic,eepic}
\usepackage{amscd,amsfonts,amsmath,amsxtra,amssymb}
\usepackage[all]{xy}
\input xypic

\newcommand{\cal}{\mathcal} 				
\renewcommand{\proof}{{\em {Proof.}\quad}}   

\newfont{\gothic}{eufm10 scaled 1100}
\newfont{\mediumgothic}{eufm10 scaled 1000}
\newfont{\smallgothic}{eufm10 scaled 900}
\newfont{\verysmallgothic}{eufm10 scaled 700}

\def \PP {{\mathbb P}}

\def \ZZ {{\mathbb Z}}

\def \cald {{\cal D}}
\def \calm {{\cal M}}
\newcommand{\mgbar}{{\overline{M}_g}}

\newcommand{\cmgbar}{{\overline{{\cal M}}_g}}

\newcommand{\cmgbarquot}{{\overline{{\cal M}}_{g^+,\mgamma}}}

\newcommand{\mgprimequot}{{M_{g^+,\mgamma}'}}

\newcommand{\cmgprimequot}{{{\cal M}_{g^+,\mgamma}'}}

\newcommand{\dpo}{{D({C_0;P_1})}}

\newcommand{\cdpo}{{{\cal D}({C_0;P_1})}}

\newcommand{\mgprime}{{M_{g^+,m}'}}

\newcommand{\hgprime}{{H_{g^+,n,m}'}}

\newcommand{\hgprimequot}{{H_{g^+,n,\mgamma}'}}

\newcommand{\hgbarquot}{{\overline{H}_ {g^+,n,\mgamma}}}

\newcommand{\cmgprime}{{{\cal M}_{g^+,m}'}}

\newcommand{\dpsplit}{{D^\times({C_0;P_1})}}

\newcommand{\thx}{{\Theta^\times}}

\newcommand{\hx}{{H_{g^+,n,m}^\times}}
\newcommand{\hxgamma}{{H_{g^+,n,\mgamma}^\times}}
\newcommand{\mxgamma}{{M_{g^+,\mgamma}^\times}}

\newcommand{\hgammax}{{H_{g^+,n,\mgamma}^\times}}

\newcommand{\mx}{{M_{g^+,m}^\times}}

\newcommand{\mgammax}{{M_{g^+,\mgamma}^\times}}

\newcommand{\kpo}{{K({C_0;P_1})}}
\newcommand{\kx}{{K^\times(C_0;P_1)}}

\newcommand{\kpx}{{K^\times}}

\newcommand{\dx}{{D^\times(C_0;P_1)}}

\newcommand{\dpx}{{D^\times}}

\newcommand{\cdx}{{{\cal D}^\times(C_0;P_1)}}

\newcommand{\cmgammax}{{{\cal M}_{g^+,\mgamma}^\times}}
\newcommand{\cmgx}{{{\cal M}_{g^+,m}^\times}}

\newcommand{\cmgammadetailx}{{{\cal M}{}^\times_{g^+,n,m/\Gamma(C_0^-;{\cal Q})}}}

\newcommand{\gammaco}{\Gamma(C_0^-;\calq)}

\newcommand{\ecs}{{E_{C/S}}}
\newcommand{\pcs}{{\phi_{C/S}}}
\newcommand{\prcs}{{p_{_{C/S}}}}

\newcommand{\cmsbar}{{\overline{\cal M}}}
\newcommand{\mgamma}{{m/\Gamma}}

\newcommand{\calq}{{\cal Q}}
\newcommand{\msbar}{\overline{{M}}}				
\newcommand{\mgstack}[1]{{\overline{\cal M}_{{#1}} }}		

\newcommand{\FDPtilde}{{\tilde{\cald}(C_0;P_1)}}	
\newcommand{\FDPopen}{{\cdpo}}					
	
\newcommand{\FDPcopen}{{{\cald}{}^\times(C_0;P_1)}}	
\newcommand{\hilb}[1]{\overline{H}_{{{#1}}}}
\newcommand{\CQ}{{\cal Q}}

\renewcommand{\epsilon}{\varepsilon}
\renewcommand{\rho}{\varrho}
\renewcommand{\theta}{\vartheta}

\newcommand{\Hilb}{\text{\rm Hilb}\,}

\newcommand{\Aut}{\text{\rm Aut}  }
\newcommand{\Stab}{\text{\rm Stab}}
\newcommand{\id}{\text{\rm id}}
\newcommand{\pr}{\text{\rm pr}}

\newcommand{\PGL}{\text{\rm PGL}}
\newcommand{\Spec}{\text{\rm Spec\,}}

\newcommand{\ebew}{\hfill$\Box$ \par}
\renewcommand{\Text}[1]{\mbox{\rm #1}}

\newcommand{\proofof}[1]{{\em {Proof of #1.}\quad}}

\newtheorem{thm}{Theorem}[section]
\newtheorem{defi}[thm]{Definition}
\newtheorem{prop}[thm]{Proposition}
\newtheorem{ex}[thm]{Example}
\newtheorem{rk}[thm]{Remark}

\newtheorem{cor}[thm]{Corollary}
\newtheorem{lemma}[thm]{Lemma}
\newtheorem{notation}[thm]{Notation}
\newtheorem{abschnitt}[thm]{}

\renewcommand{\footnote}[1]{{}}
\begin{document}

\title{The one-dimensional stratum in the boundary of the moduli stack of stable curves}

\author{J\"org Zintl}
\date{21.12.2007}
\thanks{\hspace{-4mm}{\bf 2000 Mathematics Subject Classification:} 14H10 (Primary), 14H37 (Secondary).\\   
{\bf Keywords:} Moduli of curves, moduli stacks, automorphisms of pointed stable curves}
\maketitle

\markright{The One-dimensional Stratum In The Boundary Of The Moduli Stack Of Stable Curves}

\begin{abstract}
It is well-known that the moduli space $\mgbar$ of  Deligne-Mumford  stable curves of genus $g$ admits a stratification by the loci of stable curves with a fixed number $i$ of nodes, where $0\le i \le 3g-3$. There is an analogous stratification of the associated moduli stack $\cmgbar$. 

In this paper we are interested in that particular  stratum of the moduli stack, which corresponds to stable curves with exactly $3g-4$ nodes. The irreducible components of this stratum are one-dimensional substacks of $\cmgbar$. We show how these substacks can be related to simpler moduli stacks of (permutation classes of) pointed stable curves. Furthermore, we use this to construct all of the components of this boundary stratum generically in a new way as explicit quotient stacks.   
\end{abstract}

\section{Introduction}

A stable curve $f: C \rightarrow \Spec(k)$ of genus $g\ge 2$ in the sense of Deligne and Mumford \cite{DM} may have at most $3g-3$ nodes. The locus of curves with exactly $i$ nodes is a subscheme $M_g^{(i)}$ of the moduli space $\mgbar$, and it is of dimension $3g-3-i$. Thus there is a stratification of the moduli space of Deligne-Mumford stable curves 
\[ \mgbar =  M_g^{(0)} \cup \ldots \cup M_g^{(3g-3)} . \]
If the maximal number of nodes is attained by a stable curve, then this curve decomposes into $2g-2$ irreducible components. Up to isomorphisms, there are only finitely many such curves. Therefore their locus $M_g^{(3g-3)}$ is discrete in $\mgbar$. 

Our objects of study in this paper are the connected components of the subscheme   ${M}_g^{(3g-4)}$. Each connected component $D$ of this stratum can be identified with a corresponding much simpler moduli space $M$. Before we make this formal in section \ref{sectauto}, let us briefly sketch the idea behind this approach. 

Let $f_0: C_0 \rightarrow \Spec(k)$ be a fixed stable curve of genus $g$ with $3g-3$ nodes, and let one node $P_1$ of $C_0$ be fixed. The pair $(C_0,P_1)$ distinguishes one connected component $D$ of the stratum ${M}_g^{(3g-4)}$, as we will explain in section \ref{sectmodprob}. Indeed, for any deformation of $C_0$, which preserves all nodes of $C_0$ except $P_1$, all fibres of the deformation are represented by points of $D$. Each connected component of ${M}_g^{(3g-4)}$ can be obtained in this way. 

It is known that the respective closures of these one-dimensional subschemes $D$ in $\mgbar$ are nonsingular rational curves, meeting each other transversally exactly in the points of $M_g^{(3g-3)}$. In particular, the connected components $D$ are simultaneously the irreducible components of $M_g^{(3g-4)}$.  A closed point of $D$ represents a stable curve $f: C\rightarrow \Spec(k)$ with exactly $3g-4$ nodes, and at least $2g-3$ irreducible components. Each of the irreducible components of $C$ can be considered as a pointed stable curve itself, the special points being the points of intersection with other components. There is exactly one irreducible component $C^+$ of $C$, which has not the maximal number of nodes, when considered as a pointed stable curve. The closure of the complement of $C^+$ in $C$ is a pointed prestable curve $C_0^-$, which is the same for all curves represented by closed points of $D$. In fact, $C_0^-$ is the closure of  the complement of $C_0^+$ in $C_0$, where $C_0^+$ is the union of all irreducible components of $C_0$ which contain the node $P_1$. 

\begin{ex}\em
Consider the stable curve $C_0$ of genus $g=3$ as pictured below. It attains the  maximal number of nodes $3g-3 =6$. Each component of $C_0$ is a rational curve. Let $P_1, \ldots, P_6$ be an enumeration of the nodes of $C_0$ as indicated. The node $P_1$ determines a decomposition of $C_0$ into two $2$-pointed subcurves $C_0^+$ and $C_0^+$ as shown on the right hand side of the arrow in the picture. 
 
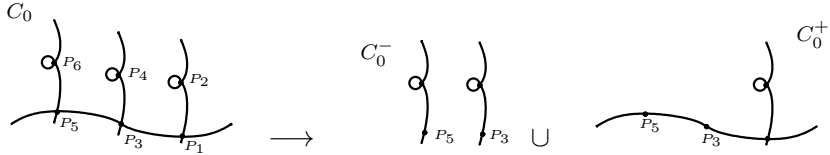
\begin{tabular}{ccccc}
\setlength{\unitlength}{0.0001in}
\begingroup\makeatletter\ifx\SetFigFont\undefined%
\gdef\SetFigFont#1#2#3#4#5{%
  \reset@font\fontsize{#1}{#2pt}%
  \fontfamily{#3}\fontseries{#4}\fontshape{#5}%
  \selectfont}%
\fi\endgroup%
{\renewcommand{\dashlinestretch}{30}
\begin{picture}(10866,6198)(0,-10)
\thicklines
\put(2283,1650){\blacken\ellipse{150}{150}}
\put(2283,1650){\ellipse{150}{150}}
\put(2133,4050){\blacken\ellipse{150}{150}}
\put(2133,4050){\ellipse{150}{150}}
\path(2133,4125)(2133,4124)(2134,4123)
	(2136,4120)(2139,4115)(2143,4109)
	(2149,4100)(2156,4088)(2164,4075)
	(2173,4059)(2184,4040)(2196,4019)
	(2208,3995)(2222,3970)(2236,3942)
	(2250,3913)(2264,3881)(2279,3848)
	(2293,3813)(2307,3776)(2321,3736)
	(2335,3695)(2348,3651)(2360,3605)
	(2372,3555)(2383,3502)(2393,3445)
	(2403,3385)(2411,3319)(2419,3249)
	(2425,3175)(2429,3096)(2432,3012)
	(2433,2925)(2432,2853)(2430,2781)
	(2427,2709)(2423,2637)(2418,2567)
	(2412,2499)(2406,2432)(2398,2367)
	(2390,2303)(2382,2241)(2373,2180)
	(2364,2120)(2354,2062)(2343,2005)
	(2333,1948)(2322,1893)(2311,1839)
	(2299,1785)(2288,1733)(2276,1682)
	(2264,1631)(2253,1582)(2241,1535)
	(2230,1489)(2219,1445)(2208,1403)
	(2198,1364)(2188,1327)(2179,1293)
	(2171,1262)(2163,1234)(2157,1210)
	(2151,1189)(2146,1171)(2142,1157)
	(2139,1146)(2136,1137)(2135,1131)
	(2134,1128)(2133,1126)(2133,1125)
\path(2133,4050)(2134,4051)(2135,4053)
	(2138,4056)(2142,4062)(2148,4069)
	(2155,4079)(2165,4092)(2175,4108)
	(2188,4125)(2202,4146)(2216,4168)
	(2232,4192)(2248,4219)(2264,4247)
	(2281,4277)(2297,4308)(2313,4341)
	(2328,4376)(2343,4413)(2358,4453)
	(2371,4494)(2384,4539)(2395,4586)
	(2406,4638)(2415,4693)(2422,4752)
	(2428,4815)(2432,4881)(2433,4950)
	(2432,5012)(2429,5074)(2425,5135)
	(2419,5195)(2411,5252)(2403,5308)
	(2393,5362)(2383,5414)(2372,5464)
	(2360,5513)(2348,5561)(2335,5607)
	(2321,5652)(2307,5696)(2293,5739)
	(2279,5781)(2264,5822)(2250,5861)
	(2236,5899)(2222,5935)(2208,5969)
	(2196,6000)(2184,6029)(2173,6055)
	(2164,6078)(2156,6097)(2149,6113)
	(2143,6126)(2139,6136)(2136,6142)
	(2134,6147)(2133,6149)(2133,6150)
\put(5283,3450){\blacken\ellipse{150}{150}}
\put(5283,3450){\ellipse{150}{150}}
\put(5433,1050){\blacken\ellipse{150}{150}}
\put(5433,1050){\ellipse{150}{150}}
\put(8358,3075){\blacken\ellipse{150}{150}}
\put(8358,3075){\ellipse{150}{150}}
\put(8433,450){\blacken\ellipse{150}{150}}
\put(8433,450){\ellipse{150}{150}}
\put(1833,4088){\ellipse{604}{604}}
\put(4983,3488){\ellipse{604}{604}}
\put(8058,3113){\ellipse{604}{604}}
\path(33,1050)(34,1051)(36,1052)
	(39,1054)(43,1057)(49,1061)
	(57,1066)(66,1072)(77,1080)
	(91,1089)(106,1099)(125,1110)
	(145,1123)(168,1137)(193,1152)
	(220,1168)(250,1185)(282,1204)
	(316,1223)(353,1242)(391,1263)
	(432,1283)(475,1304)(519,1326)
	(566,1347)(614,1369)(665,1390)
	(717,1411)(771,1432)(827,1453)
	(886,1473)(946,1492)(1009,1511)
	(1074,1530)(1142,1547)(1213,1564)
	(1286,1580)(1363,1595)(1443,1609)
	(1527,1622)(1615,1634)(1707,1644)
	(1803,1653)(1904,1661)(2010,1668)
	(2121,1672)(2237,1675)(2358,1676)
	(2484,1675)(2615,1673)(2750,1667)
	(2890,1660)(3033,1650)(3145,1641)
	(3258,1630)(3370,1618)(3482,1605)
	(3592,1591)(3700,1577)(3806,1561)
	(3909,1545)(4008,1528)(4104,1511)
	(4197,1493)(4286,1475)(4370,1457)
	(4451,1438)(4528,1420)(4601,1401)
	(4670,1382)(4735,1364)(4796,1345)
	(4854,1326)(4909,1307)(4960,1288)
	(5009,1270)(5054,1251)(5097,1233)
	(5138,1214)(5176,1196)(5212,1177)
	(5247,1159)(5280,1141)(5312,1122)
	(5343,1104)(5374,1086)(5403,1068)
	(5433,1050)(5463,1032)(5492,1014)
	(5523,996)(5554,978)(5586,959)
	(5619,941)(5654,923)(5690,904)
	(5728,886)(5769,867)(5812,849)
	(5857,830)(5906,812)(5957,793)
	(6012,774)(6070,755)(6131,736)
	(6196,718)(6265,699)(6338,680)
	(6415,662)(6496,643)(6580,625)
	(6669,607)(6762,589)(6858,572)
	(6957,555)(7060,539)(7166,523)
	(7274,509)(7384,495)(7496,482)
	(7608,470)(7721,459)(7833,450)
	(7976,440)(8116,433)(8251,427)
	(8382,425)(8508,424)(8629,425)
	(8745,428)(8856,432)(8962,439)
	(9063,447)(9159,456)(9251,466)
	(9339,478)(9423,491)(9503,505)
	(9580,520)(9653,536)(9724,553)
	(9792,570)(9857,589)(9920,608)
	(9980,627)(10039,647)(10095,668)
	(10149,689)(10201,710)(10252,731)
	(10300,753)(10347,774)(10391,796)
	(10434,817)(10475,837)(10513,858)
	(10550,877)(10584,896)(10616,915)
	(10646,932)(10673,948)(10698,963)
	(10721,977)(10741,990)(10760,1001)
	(10775,1011)(10789,1020)(10800,1028)
	(10809,1034)(10817,1039)(10823,1043)
	(10827,1046)(10830,1048)(10832,1049)(10833,1050)
\path(5283,3525)(5283,3524)(5284,3523)
	(5286,3520)(5289,3515)(5293,3509)
	(5299,3500)(5306,3488)(5314,3475)
	(5323,3459)(5334,3440)(5346,3419)
	(5358,3395)(5372,3370)(5386,3342)
	(5400,3313)(5414,3281)(5429,3248)
	(5443,3213)(5457,3176)(5471,3136)
	(5485,3095)(5498,3051)(5510,3005)
	(5522,2955)(5533,2902)(5543,2845)
	(5553,2785)(5561,2719)(5569,2649)
	(5575,2575)(5579,2496)(5582,2412)
	(5583,2325)(5582,2253)(5580,2181)
	(5577,2109)(5573,2037)(5568,1967)
	(5562,1899)(5556,1832)(5548,1767)
	(5540,1703)(5532,1641)(5523,1580)
	(5514,1520)(5504,1462)(5493,1405)
	(5483,1348)(5472,1293)(5461,1239)
	(5449,1185)(5438,1133)(5426,1082)
	(5414,1031)(5403,982)(5391,935)
	(5380,889)(5369,845)(5358,803)
	(5348,764)(5338,727)(5329,693)
	(5321,662)(5313,634)(5307,610)
	(5301,589)(5296,571)(5292,557)
	(5289,546)(5286,537)(5285,531)
	(5284,528)(5283,526)(5283,525)
\path(5283,3450)(5284,3451)(5285,3453)
	(5288,3456)(5292,3462)(5298,3469)
	(5305,3479)(5315,3492)(5325,3508)
	(5338,3525)(5352,3546)(5366,3568)
	(5382,3592)(5398,3619)(5414,3647)
	(5431,3677)(5447,3708)(5463,3741)
	(5478,3776)(5493,3813)(5508,3853)
	(5521,3894)(5534,3939)(5545,3986)
	(5556,4038)(5565,4093)(5572,4152)
	(5578,4215)(5582,4281)(5583,4350)
	(5582,4412)(5579,4474)(5575,4535)
	(5569,4595)(5561,4652)(5553,4708)
	(5543,4762)(5533,4814)(5522,4864)
	(5510,4913)(5498,4961)(5485,5007)
	(5471,5052)(5457,5096)(5443,5139)
	(5429,5181)(5414,5222)(5400,5261)
	(5386,5299)(5372,5335)(5358,5369)
	(5346,5400)(5334,5429)(5323,5455)
	(5314,5478)(5306,5497)(5299,5513)
	(5293,5526)(5289,5536)(5286,5542)
	(5284,5547)(5283,5549)(5283,5550)
\path(8358,3150)(8358,3149)(8359,3148)
	(8361,3145)(8364,3140)(8368,3134)
	(8374,3125)(8381,3113)(8389,3100)
	(8398,3084)(8409,3065)(8421,3044)
	(8433,3020)(8447,2995)(8461,2967)
	(8475,2938)(8489,2906)(8504,2873)
	(8518,2838)(8532,2801)(8546,2761)
	(8560,2720)(8573,2676)(8585,2630)
	(8597,2580)(8608,2527)(8618,2470)
	(8628,2410)(8636,2344)(8644,2274)
	(8650,2200)(8654,2121)(8657,2037)
	(8658,1950)(8657,1878)(8655,1806)
	(8652,1734)(8648,1662)(8643,1592)
	(8637,1524)(8631,1457)(8623,1392)
	(8615,1328)(8607,1266)(8598,1205)
	(8589,1145)(8579,1087)(8568,1030)
	(8558,973)(8547,918)(8536,864)
	(8524,810)(8513,758)(8501,707)
	(8489,656)(8478,607)(8466,560)
	(8455,514)(8444,470)(8433,428)
	(8423,389)(8413,352)(8404,318)
	(8396,287)(8388,259)(8382,235)
	(8376,214)(8371,196)(8367,182)
	(8364,171)(8361,162)(8360,156)
	(8359,153)(8358,151)(8358,150)
\path(8358,3075)(8359,3076)(8360,3078)
	(8363,3081)(8367,3087)(8373,3094)
	(8380,3104)(8390,3117)(8400,3133)
	(8413,3150)(8427,3171)(8441,3193)
	(8457,3217)(8473,3244)(8489,3272)
	(8506,3302)(8522,3333)(8538,3366)
	(8553,3401)(8568,3438)(8583,3478)
	(8596,3519)(8609,3564)(8620,3611)
	(8631,3663)(8640,3718)(8647,3777)
	(8653,3840)(8657,3906)(8658,3975)
	(8657,4037)(8654,4099)(8650,4160)
	(8644,4220)(8636,4277)(8628,4333)
	(8618,4387)(8608,4439)(8597,4489)
	(8585,4538)(8573,4586)(8560,4632)
	(8546,4677)(8532,4721)(8518,4764)
	(8504,4806)(8489,4847)(8475,4886)
	(8461,4924)(8447,4960)(8433,4994)
	(8421,5025)(8409,5054)(8398,5080)
	(8389,5103)(8381,5122)(8374,5138)
	(8368,5151)(8364,5161)(8361,5167)
	(8359,5172)(8358,5174)(8358,5175)
\put(8733,3075){\makebox(0,0)[lb]{\smash{{{\SetFigFont{5}{6.0}{\rmdefault}{\mddefault}{\updefault}$P_2$}}}}}
\put(8508,-200){\makebox(0,0)[lb]{\smash{{{\SetFigFont{5}{6.0}{\rmdefault}{\mddefault}{\updefault}$P_1$}}}}}
\put(5720,3375){\makebox(0,0)[lb]{\smash{{{\SetFigFont{5}{6.0}{\rmdefault}{\mddefault}{\updefault}$P_4$}}}}}
\put(2483,3900){\makebox(0,0)[lb]{\smash{{{\SetFigFont{5}{6.0}{\rmdefault}{\mddefault}{\updefault}$P_6$}}}}}
\put(2358,900){\makebox(0,0)[lb]{\smash{{{\SetFigFont{5}{6.0}{\rmdefault}{\mddefault}{\updefault}$P_5$}}}}}
\put(5508,50){\makebox(0,0)[lb]{\smash{{{\SetFigFont{5}{6.0}{\rmdefault}{\mddefault}{\updefault}$P_3$}}}}}
\put(-200,6233){\makebox(0,0)[lb]{\smash{{{\SetFigFont{8}{9.6}{\rmdefault}{\mddefault}{\updefault}$C_0$}}}}}
\end{picture}
}
& $\longrightarrow$ &
\setlength{\unitlength}{0.0001in}
\begingroup\makeatletter\ifx\SetFigFont\undefined%
\gdef\SetFigFont#1#2#3#4#5{%
  \reset@font\fontsize{#1}{#2pt}%
  \fontfamily{#3}\fontseries{#4}\fontshape{#5}%
  \selectfont}%
\fi\endgroup%
{\renewcommand{\dashlinestretch}{30}
\begin{picture}(6592,5181)(0,-10)
\thicklines
\put(3150,633){\blacken\ellipse{150}{150}}
\put(3150,633){\ellipse{150}{150}}
\put(3000,3033){\blacken\ellipse{150}{150}}
\put(3000,3033){\ellipse{150}{150}}
\path(3000,3108)(3000,3107)(3001,3106)
	(3003,3103)(3006,3098)(3010,3092)
	(3016,3083)(3023,3071)(3031,3058)
	(3040,3042)(3051,3023)(3063,3002)
	(3075,2978)(3089,2953)(3103,2925)
	(3117,2896)(3131,2864)(3146,2831)
	(3160,2796)(3174,2759)(3188,2719)
	(3202,2678)(3215,2634)(3227,2588)
	(3239,2538)(3250,2485)(3260,2428)
	(3270,2368)(3278,2302)(3286,2232)
	(3292,2158)(3296,2079)(3299,1995)
	(3300,1908)(3299,1836)(3297,1764)
	(3294,1692)(3290,1620)(3285,1550)
	(3279,1482)(3273,1415)(3265,1350)
	(3257,1286)(3249,1224)(3240,1163)
	(3231,1103)(3221,1045)(3210,988)
	(3200,931)(3189,876)(3178,822)
	(3166,768)(3155,716)(3143,665)
	(3131,614)(3120,565)(3108,518)
	(3097,472)(3086,428)(3075,386)
	(3065,347)(3055,310)(3046,276)
	(3038,245)(3030,217)(3024,193)
	(3018,172)(3013,154)(3009,140)
	(3006,129)(3003,120)(3002,114)
	(3001,111)(3000,109)(3000,108)
\path(3000,3033)(3001,3034)(3002,3036)
	(3005,3039)(3009,3045)(3015,3052)
	(3022,3062)(3032,3075)(3042,3091)
	(3055,3108)(3069,3129)(3083,3151)
	(3099,3175)(3115,3202)(3131,3230)
	(3148,3260)(3164,3291)(3180,3324)
	(3195,3359)(3210,3396)(3225,3436)
	(3238,3477)(3251,3522)(3262,3569)
	(3273,3621)(3282,3676)(3289,3735)
	(3295,3798)(3299,3864)(3300,3933)
	(3299,3995)(3296,4057)(3292,4118)
	(3286,4178)(3278,4235)(3270,4291)
	(3260,4345)(3250,4397)(3239,4447)
	(3227,4496)(3215,4544)(3202,4590)
	(3188,4635)(3174,4679)(3160,4722)
	(3146,4764)(3131,4805)(3117,4844)
	(3103,4882)(3089,4918)(3075,4952)
	(3063,4983)(3051,5012)(3040,5038)
	(3031,5061)(3023,5080)(3016,5096)
	(3010,5109)(3006,5119)(3003,5125)
	(3001,5130)(3000,5132)(3000,5133)
\put(5850,2958){\blacken\ellipse{150}{150}}
\put(5850,2958){\ellipse{150}{150}}
\put(6000,558){\blacken\ellipse{150}{150}}
\put(6000,558){\ellipse{150}{150}}
\path(5850,3033)(5850,3032)(5851,3031)
	(5853,3028)(5856,3023)(5860,3017)
	(5866,3008)(5873,2996)(5881,2983)
	(5890,2967)(5901,2948)(5913,2927)
	(5925,2903)(5939,2878)(5953,2850)
	(5967,2821)(5981,2789)(5996,2756)
	(6010,2721)(6024,2684)(6038,2644)
	(6052,2603)(6065,2559)(6077,2513)
	(6089,2463)(6100,2410)(6110,2353)
	(6120,2293)(6128,2227)(6136,2157)
	(6142,2083)(6146,2004)(6149,1920)
	(6150,1833)(6149,1761)(6147,1689)
	(6144,1617)(6140,1545)(6135,1475)
	(6129,1407)(6123,1340)(6115,1275)
	(6107,1211)(6099,1149)(6090,1088)
	(6081,1028)(6071,970)(6060,913)
	(6050,856)(6039,801)(6028,747)
	(6016,693)(6005,641)(5993,590)
	(5981,539)(5970,490)(5958,443)
	(5947,397)(5936,353)(5925,311)
	(5915,272)(5905,235)(5896,201)
	(5888,170)(5880,142)(5874,118)
	(5868,97)(5863,79)(5859,65)
	(5856,54)(5853,45)(5852,39)
	(5851,36)(5850,34)(5850,33)
\path(5850,2958)(5851,2959)(5852,2961)
	(5855,2964)(5859,2970)(5865,2977)
	(5872,2987)(5882,3000)(5892,3016)
	(5905,3033)(5919,3054)(5933,3076)
	(5949,3100)(5965,3127)(5981,3155)
	(5998,3185)(6014,3216)(6030,3249)
	(6045,3284)(6060,3321)(6075,3361)
	(6088,3402)(6101,3447)(6112,3494)
	(6123,3546)(6132,3601)(6139,3660)
	(6145,3723)(6149,3789)(6150,3858)
	(6149,3920)(6146,3982)(6142,4043)
	(6136,4103)(6128,4160)(6120,4216)
	(6110,4270)(6100,4322)(6089,4372)
	(6077,4421)(6065,4469)(6052,4515)
	(6038,4560)(6024,4604)(6010,4647)
	(5996,4689)(5981,4730)(5967,4769)
	(5953,4807)(5939,4843)(5925,4877)
	(5913,4908)(5901,4937)(5890,4963)
	(5881,4986)(5873,5005)(5866,5021)
	(5860,5034)(5856,5044)(5853,5050)
	(5851,5055)(5850,5057)(5850,5058)
\put(2700,3071){\ellipse{604}{604}}
\put(5550,2996){\ellipse{604}{604}}
\put(3600,333){\makebox(0,0)[lb]{\smash{{{\SetFigFont{5}{6.0}{\rmdefault}{\mddefault}{\updefault}$P_5$}}}}}
\put(6300,333){\makebox(0,0)[lb]{\smash{{{\SetFigFont{5}{6.0}{\rmdefault}{\mddefault}{\updefault}$P_3$}}}}}
\put(0,4233){\makebox(0,0)[lb]{\smash{{{\SetFigFont{8}{9.6}{\rmdefault}{\mddefault}{\updefault}$C_0^{-}$}}}}}
\end{picture}
}
& $\cup$  &
\setlength{\unitlength}{0.0001in}
\begingroup\makeatletter\ifx\SetFigFont\undefined%
\gdef\SetFigFont#1#2#3#4#5{%
  \reset@font\fontsize{#1}{#2pt}%
  \fontfamily{#3}\fontseries{#4}\fontshape{#5}%
  \selectfont}%
\fi\endgroup%
{\renewcommand{\dashlinestretch}{30}
\begin{picture}(12180,5658)(0,-10)
\thicklines
\put(5433,933){\blacken\ellipse{150}{150}}
\put(5433,933){\ellipse{150}{150}}
\put(8358,2958){\blacken\ellipse{150}{150}}
\put(8358,2958){\ellipse{150}{150}}
\put(8433,333){\blacken\ellipse{150}{150}}
\put(8433,333){\ellipse{150}{150}}
\put(8058,2996){\ellipse{604}{604}}
\put(2433,1533){\blacken\ellipse{150}{150}}
\put(2433,1533){\ellipse{150}{150}}
\path(33,933)(34,934)(36,935)
	(39,937)(43,940)(49,944)
	(57,949)(66,955)(77,963)
	(91,972)(106,982)(125,993)
	(145,1006)(168,1020)(193,1035)
	(220,1051)(250,1068)(282,1087)
	(316,1106)(353,1125)(391,1146)
	(432,1166)(475,1187)(519,1209)
	(566,1230)(614,1252)(665,1273)
	(717,1294)(771,1315)(827,1336)
	(886,1356)(946,1375)(1009,1394)
	(1074,1413)(1142,1430)(1213,1447)
	(1286,1463)(1363,1478)(1443,1492)
	(1527,1505)(1615,1517)(1707,1527)
	(1803,1536)(1904,1544)(2010,1551)
	(2121,1555)(2237,1558)(2358,1559)
	(2484,1558)(2615,1556)(2750,1550)
	(2890,1543)(3033,1533)(3145,1524)
	(3258,1513)(3370,1501)(3482,1488)
	(3592,1474)(3700,1460)(3806,1444)
	(3909,1428)(4008,1411)(4104,1394)
	(4197,1376)(4286,1358)(4370,1340)
	(4451,1321)(4528,1303)(4601,1284)
	(4670,1265)(4735,1247)(4796,1228)
	(4854,1209)(4909,1190)(4960,1171)
	(5009,1153)(5054,1134)(5097,1116)
	(5138,1097)(5176,1079)(5212,1060)
	(5247,1042)(5280,1024)(5312,1005)
	(5343,987)(5374,969)(5403,951)
	(5433,933)(5463,915)(5492,897)
	(5523,879)(5554,861)(5586,842)
	(5619,824)(5654,806)(5690,787)
	(5728,769)(5769,750)(5812,732)
	(5857,713)(5906,695)(5957,676)
	(6012,657)(6070,638)(6131,619)
	(6196,601)(6265,582)(6338,563)
	(6415,545)(6496,526)(6580,508)
	(6669,490)(6762,472)(6858,455)
	(6957,438)(7060,422)(7166,406)
	(7274,392)(7384,378)(7496,365)
	(7608,353)(7721,342)(7833,333)
	(7976,323)(8116,316)(8251,310)
	(8382,308)(8508,307)(8629,308)
	(8745,311)(8856,315)(8962,322)
	(9063,330)(9159,339)(9251,349)
	(9339,361)(9423,374)(9503,388)
	(9580,403)(9653,419)(9724,436)
	(9792,453)(9857,472)(9920,491)
	(9980,510)(10039,530)(10095,551)
	(10149,572)(10201,593)(10252,614)
	(10300,636)(10347,657)(10391,679)
	(10434,700)(10475,720)(10513,741)
	(10550,760)(10584,779)(10616,798)
	(10646,815)(10673,831)(10698,846)
	(10721,860)(10741,873)(10760,884)
	(10775,894)(10789,903)(10800,911)
	(10809,917)(10817,922)(10823,926)
	(10827,929)(10830,931)(10832,932)(10833,933)
\path(8358,3033)(8358,3032)(8359,3031)
	(8361,3028)(8364,3023)(8368,3017)
	(8374,3008)(8381,2996)(8389,2983)
	(8398,2967)(8409,2948)(8421,2927)
	(8433,2903)(8447,2878)(8461,2850)
	(8475,2821)(8489,2789)(8504,2756)
	(8518,2721)(8532,2684)(8546,2644)
	(8560,2603)(8573,2559)(8585,2513)
	(8597,2463)(8608,2410)(8618,2353)
	(8628,2293)(8636,2227)(8644,2157)
	(8650,2083)(8654,2004)(8657,1920)
	(8658,1833)(8657,1761)(8655,1689)
	(8652,1617)(8648,1545)(8643,1475)
	(8637,1407)(8631,1340)(8623,1275)
	(8615,1211)(8607,1149)(8598,1088)
	(8589,1028)(8579,970)(8568,913)
	(8558,856)(8547,801)(8536,747)
	(8524,693)(8513,641)(8501,590)
	(8489,539)(8478,490)(8466,443)
	(8455,397)(8444,353)(8433,311)
	(8423,272)(8413,235)(8404,201)
	(8396,170)(8388,142)(8382,118)
	(8376,97)(8371,79)(8367,65)
	(8364,54)(8361,45)(8360,39)
	(8359,36)(8358,34)(8358,33)
\path(8358,2958)(8359,2959)(8360,2961)
	(8363,2964)(8367,2970)(8373,2977)
	(8380,2987)(8390,3000)(8400,3016)
	(8413,3033)(8427,3054)(8441,3076)
	(8457,3100)(8473,3127)(8489,3155)
	(8506,3185)(8522,3216)(8538,3249)
	(8553,3284)(8568,3321)(8583,3361)
	(8596,3402)(8609,3447)(8620,3494)
	(8631,3546)(8640,3601)(8647,3660)
	(8653,3723)(8657,3789)(8658,3858)
	(8657,3920)(8654,3982)(8650,4043)
	(8644,4103)(8636,4160)(8628,4216)
	(8618,4270)(8608,4322)(8597,4372)
	(8585,4421)(8573,4469)(8560,4515)
	(8546,4560)(8532,4604)(8518,4647)
	(8504,4689)(8489,4730)(8475,4769)
	(8461,4807)(8447,4843)(8433,4877)
	(8421,4908)(8409,4937)(8398,4963)
	(8389,4986)(8381,5005)(8374,5021)
	(8368,5034)(8364,5044)(8361,5050)
	(8359,5055)(8358,5057)(8358,5058)
\put(5133,133){\makebox(0,0)[lb]{\smash{{{\SetFigFont{5}{6.0}{\rmdefault}{\mddefault}{\updefault}$P_3$}}}}}
\put(2133,833){\makebox(0,0)[lb]{\smash{{{\SetFigFont{5}{6.0}{\rmdefault}{\mddefault}{\updefault}$P_5$}}}}}
\put(9858,5133){\makebox(0,0)[lb]{\smash{{{\SetFigFont{8}{9.6}{\rmdefault}{\mddefault}{\updefault}$C_0^{+}$}}}}}
\put(8700,-300){\makebox(0,0)[lb]{\smash{{{\SetFigFont{5}{6.0}{\rmdefault}{\mddefault}{\updefault}$$}}}}}
\end{picture}
}
\end{tabular}

A curve $C$, obtained by a 
deformation of $C_0$ preserving all nodes except $P_1$, looks schematically like this.

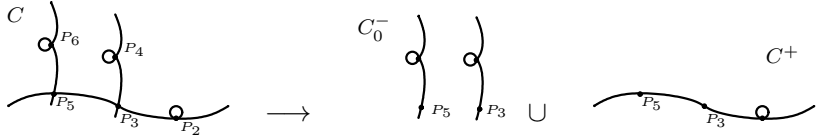
\begin{tabular}{ccccc}
\setlength{\unitlength}{0.0001in}
\begingroup\makeatletter\ifx\SetFigFont\undefined%
\gdef\SetFigFont#1#2#3#4#5{%
  \reset@font\fontsize{#1}{#2pt}%
  \fontfamily{#3}\fontseries{#4}\fontshape{#5}%
  \selectfont}%
\fi\endgroup%
{\renewcommand{\dashlinestretch}{30}
\begin{picture}(10908,5846)(0,-10)
\thicklines
\put(2325,1298){\blacken\ellipse{150}{150}}
\put(2325,1298){\ellipse{150}{150}}
\put(2175,3698){\blacken\ellipse{150}{150}}
\put(2175,3698){\ellipse{150}{150}}
\path(2175,3773)(2175,3772)(2176,3771)
	(2178,3768)(2181,3763)(2185,3757)
	(2191,3748)(2198,3736)(2206,3723)
	(2215,3707)(2226,3688)(2238,3667)
	(2250,3643)(2264,3618)(2278,3590)
	(2292,3561)(2306,3529)(2321,3496)
	(2335,3461)(2349,3424)(2363,3384)
	(2377,3343)(2390,3299)(2402,3253)
	(2414,3203)(2425,3150)(2435,3093)
	(2445,3033)(2453,2967)(2461,2897)
	(2467,2823)(2471,2744)(2474,2660)
	(2475,2573)(2474,2501)(2472,2429)
	(2469,2357)(2465,2285)(2460,2215)
	(2454,2147)(2448,2080)(2440,2015)
	(2432,1951)(2424,1889)(2415,1828)
	(2406,1768)(2396,1710)(2385,1653)
	(2375,1596)(2364,1541)(2353,1487)
	(2341,1433)(2330,1381)(2318,1330)
	(2306,1279)(2295,1230)(2283,1183)
	(2272,1137)(2261,1093)(2250,1051)
	(2240,1012)(2230,975)(2221,941)
	(2213,910)(2205,882)(2199,858)
	(2193,837)(2188,819)(2184,805)
	(2181,794)(2178,785)(2177,779)
	(2176,776)(2175,774)(2175,773)
\path(2175,3698)(2176,3699)(2177,3701)
	(2180,3704)(2184,3710)(2190,3717)
	(2197,3727)(2207,3740)(2217,3756)
	(2230,3773)(2244,3794)(2258,3816)
	(2274,3840)(2290,3867)(2306,3895)
	(2323,3925)(2339,3956)(2355,3989)
	(2370,4024)(2385,4061)(2400,4101)
	(2413,4142)(2426,4187)(2437,4234)
	(2448,4286)(2457,4341)(2464,4400)
	(2470,4463)(2474,4529)(2475,4598)
	(2474,4660)(2471,4722)(2467,4783)
	(2461,4843)(2453,4900)(2445,4956)
	(2435,5010)(2425,5062)(2414,5112)
	(2402,5161)(2390,5209)(2377,5255)
	(2363,5300)(2349,5344)(2335,5387)
	(2321,5429)(2306,5470)(2292,5509)
	(2278,5547)(2264,5583)(2250,5617)
	(2238,5648)(2226,5677)(2215,5703)
	(2206,5726)(2198,5745)(2191,5761)
	(2185,5774)(2181,5784)(2178,5790)
	(2176,5795)(2175,5797)(2175,5798)
\put(5325,3098){\blacken\ellipse{150}{150}}
\put(5325,3098){\ellipse{150}{150}}
\put(5475,698){\blacken\ellipse{150}{150}}
\put(5475,698){\ellipse{150}{150}}
\put(1875,3736){\ellipse{604}{604}}
\put(5025,3136){\ellipse{604}{604}}
\put(8325,398){\ellipse{600}{600}}
\put(8325,98){\blacken\ellipse{150}{150}}
\put(8325,98){\ellipse{150}{150}}
\path(75,698)(76,699)(78,700)
	(81,702)(85,705)(91,709)
	(99,714)(108,720)(119,728)
	(133,737)(148,747)(167,758)
	(187,771)(210,785)(235,800)
	(262,816)(292,833)(324,852)
	(358,871)(395,890)(433,911)
	(474,931)(517,952)(561,974)
	(608,995)(656,1017)(707,1038)
	(759,1059)(813,1080)(869,1101)
	(928,1121)(988,1140)(1051,1159)
	(1116,1178)(1184,1195)(1255,1212)
	(1328,1228)(1405,1243)(1485,1257)
	(1569,1270)(1657,1282)(1749,1292)
	(1845,1301)(1946,1309)(2052,1316)
	(2163,1320)(2279,1323)(2400,1324)
	(2526,1323)(2657,1321)(2792,1315)
	(2932,1308)(3075,1298)(3187,1289)
	(3300,1278)(3412,1266)(3524,1253)
	(3634,1239)(3742,1225)(3848,1209)
	(3951,1193)(4050,1176)(4146,1159)
	(4239,1141)(4328,1123)(4412,1105)
	(4493,1086)(4570,1068)(4643,1049)
	(4712,1030)(4777,1012)(4838,993)
	(4896,974)(4951,955)(5002,936)
	(5051,918)(5096,899)(5139,881)
	(5180,862)(5218,844)(5254,825)
	(5289,807)(5322,789)(5354,770)
	(5385,752)(5416,734)(5445,716)
	(5475,698)(5505,680)(5534,662)
	(5565,644)(5596,626)(5628,607)
	(5661,589)(5696,571)(5732,552)
	(5770,534)(5811,515)(5854,497)
	(5899,478)(5948,460)(5999,441)
	(6054,422)(6112,403)(6173,384)
	(6238,366)(6307,347)(6380,328)
	(6457,310)(6538,291)(6622,273)
	(6711,255)(6804,237)(6900,220)
	(6999,203)(7102,187)(7208,171)
	(7316,157)(7426,143)(7538,130)
	(7650,118)(7763,107)(7875,98)
	(8018,88)(8158,81)(8293,75)
	(8424,73)(8550,72)(8671,73)
	(8787,76)(8898,80)(9004,87)
	(9105,95)(9201,104)(9293,114)
	(9381,126)(9465,139)(9545,153)
	(9622,168)(9695,184)(9766,201)
	(9834,218)(9899,237)(9962,256)
	(10022,275)(10081,295)(10137,316)
	(10191,337)(10243,358)(10294,379)
	(10342,401)(10389,422)(10433,444)
	(10476,465)(10517,485)(10555,506)
	(10592,525)(10626,544)(10658,563)
	(10688,580)(10715,596)(10740,611)
	(10763,625)(10783,638)(10802,649)
	(10817,659)(10831,668)(10842,676)
	(10851,682)(10859,687)(10865,691)
	(10869,694)(10872,696)(10874,697)(10875,698)
\path(5325,3173)(5325,3172)(5326,3171)
	(5328,3168)(5331,3163)(5335,3157)
	(5341,3148)(5348,3136)(5356,3123)
	(5365,3107)(5376,3088)(5388,3067)
	(5400,3043)(5414,3018)(5428,2990)
	(5442,2961)(5456,2929)(5471,2896)
	(5485,2861)(5499,2824)(5513,2784)
	(5527,2743)(5540,2699)(5552,2653)
	(5564,2603)(5575,2550)(5585,2493)
	(5595,2433)(5603,2367)(5611,2297)
	(5617,2223)(5621,2144)(5624,2060)
	(5625,1973)(5624,1901)(5622,1829)
	(5619,1757)(5615,1685)(5610,1615)
	(5604,1547)(5598,1480)(5590,1415)
	(5582,1351)(5574,1289)(5565,1228)
	(5556,1168)(5546,1110)(5535,1053)
	(5525,996)(5514,941)(5503,887)
	(5491,833)(5480,781)(5468,730)
	(5456,679)(5445,630)(5433,583)
	(5422,537)(5411,493)(5400,451)
	(5390,412)(5380,375)(5371,341)
	(5363,310)(5355,282)(5349,258)
	(5343,237)(5338,219)(5334,205)
	(5331,194)(5328,185)(5327,179)
	(5326,176)(5325,174)(5325,173)
\path(5325,3098)(5326,3099)(5327,3101)
	(5330,3104)(5334,3110)(5340,3117)
	(5347,3127)(5357,3140)(5367,3156)
	(5380,3173)(5394,3194)(5408,3216)
	(5424,3240)(5440,3267)(5456,3295)
	(5473,3325)(5489,3356)(5505,3389)
	(5520,3424)(5535,3461)(5550,3501)
	(5563,3542)(5576,3587)(5587,3634)
	(5598,3686)(5607,3741)(5614,3800)
	(5620,3863)(5624,3929)(5625,3998)
	(5624,4060)(5621,4122)(5617,4183)
	(5611,4243)(5603,4300)(5595,4356)
	(5585,4410)(5575,4462)(5564,4512)
	(5552,4561)(5540,4609)(5527,4655)
	(5513,4700)(5499,4744)(5485,4787)
	(5471,4829)(5456,4870)(5442,4909)
	(5428,4947)(5414,4983)(5400,5017)
	(5388,5048)(5376,5077)(5365,5103)
	(5356,5126)(5348,5145)(5341,5161)
	(5335,5174)(5331,5184)(5328,5190)
	(5326,5195)(5325,5197)(5325,5198)
\put(0,4823){\makebox(0,0)[lb]{\smash{{{\SetFigFont{8}{9.6}{\rmdefault}{\mddefault}{\updefault}$C$}}}}}
\put(8508,-500){\makebox(0,0)[lb]{\smash{{{\SetFigFont{5}{6.0}{\rmdefault}{\mddefault}{\updefault}$P_2$}}}}}
\put(5720,3375){\makebox(0,0)[lb]{\smash{{{\SetFigFont{5}{6.0}{\rmdefault}{\mddefault}{\updefault}$P_4$}}}}}
\put(2583,3900){\makebox(0,0)[lb]{\smash{{{\SetFigFont{5}{6.0}{\rmdefault}{\mddefault}{\updefault}$P_6$}}}}}
\put(2358,700){\makebox(0,0)[lb]{\smash{{{\SetFigFont{5}{6.0}{\rmdefault}{\mddefault}{\updefault}$P_5$}}}}}
\put(5508,-150){\makebox(0,0)[lb]{\smash{{{\SetFigFont{5}{6.0}{\rmdefault}{\mddefault}{\updefault}$P_3$}}}}}
\end{picture}
}
& $\longrightarrow$ &
\setlength{\unitlength}{0.0001in}
\begingroup\makeatletter\ifx\SetFigFont\undefined%
\gdef\SetFigFont#1#2#3#4#5{%
  \reset@font\fontsize{#1}{#2pt}%
  \fontfamily{#3}\fontseries{#4}\fontshape{#5}%
  \selectfont}%
\fi\endgroup%
{\renewcommand{\dashlinestretch}{30}
\begin{picture}(6592,5181)(0,-10)
\thicklines
\put(3150,633){\blacken\ellipse{150}{150}}
\put(3150,633){\ellipse{150}{150}}
\put(3000,3033){\blacken\ellipse{150}{150}}
\put(3000,3033){\ellipse{150}{150}}
\path(3000,3108)(3000,3107)(3001,3106)
	(3003,3103)(3006,3098)(3010,3092)
	(3016,3083)(3023,3071)(3031,3058)
	(3040,3042)(3051,3023)(3063,3002)
	(3075,2978)(3089,2953)(3103,2925)
	(3117,2896)(3131,2864)(3146,2831)
	(3160,2796)(3174,2759)(3188,2719)
	(3202,2678)(3215,2634)(3227,2588)
	(3239,2538)(3250,2485)(3260,2428)
	(3270,2368)(3278,2302)(3286,2232)
	(3292,2158)(3296,2079)(3299,1995)
	(3300,1908)(3299,1836)(3297,1764)
	(3294,1692)(3290,1620)(3285,1550)
	(3279,1482)(3273,1415)(3265,1350)
	(3257,1286)(3249,1224)(3240,1163)
	(3231,1103)(3221,1045)(3210,988)
	(3200,931)(3189,876)(3178,822)
	(3166,768)(3155,716)(3143,665)
	(3131,614)(3120,565)(3108,518)
	(3097,472)(3086,428)(3075,386)
	(3065,347)(3055,310)(3046,276)
	(3038,245)(3030,217)(3024,193)
	(3018,172)(3013,154)(3009,140)
	(3006,129)(3003,120)(3002,114)
	(3001,111)(3000,109)(3000,108)
\path(3000,3033)(3001,3034)(3002,3036)
	(3005,3039)(3009,3045)(3015,3052)
	(3022,3062)(3032,3075)(3042,3091)
	(3055,3108)(3069,3129)(3083,3151)
	(3099,3175)(3115,3202)(3131,3230)
	(3148,3260)(3164,3291)(3180,3324)
	(3195,3359)(3210,3396)(3225,3436)
	(3238,3477)(3251,3522)(3262,3569)
	(3273,3621)(3282,3676)(3289,3735)
	(3295,3798)(3299,3864)(3300,3933)
	(3299,3995)(3296,4057)(3292,4118)
	(3286,4178)(3278,4235)(3270,4291)
	(3260,4345)(3250,4397)(3239,4447)
	(3227,4496)(3215,4544)(3202,4590)
	(3188,4635)(3174,4679)(3160,4722)
	(3146,4764)(3131,4805)(3117,4844)
	(3103,4882)(3089,4918)(3075,4952)
	(3063,4983)(3051,5012)(3040,5038)
	(3031,5061)(3023,5080)(3016,5096)
	(3010,5109)(3006,5119)(3003,5125)
	(3001,5130)(3000,5132)(3000,5133)
\put(5850,2958){\blacken\ellipse{150}{150}}
\put(5850,2958){\ellipse{150}{150}}
\put(6000,558){\blacken\ellipse{150}{150}}
\put(6000,558){\ellipse{150}{150}}
\path(5850,3033)(5850,3032)(5851,3031)
	(5853,3028)(5856,3023)(5860,3017)
	(5866,3008)(5873,2996)(5881,2983)
	(5890,2967)(5901,2948)(5913,2927)
	(5925,2903)(5939,2878)(5953,2850)
	(5967,2821)(5981,2789)(5996,2756)
	(6010,2721)(6024,2684)(6038,2644)
	(6052,2603)(6065,2559)(6077,2513)
	(6089,2463)(6100,2410)(6110,2353)
	(6120,2293)(6128,2227)(6136,2157)
	(6142,2083)(6146,2004)(6149,1920)
	(6150,1833)(6149,1761)(6147,1689)
	(6144,1617)(6140,1545)(6135,1475)
	(6129,1407)(6123,1340)(6115,1275)
	(6107,1211)(6099,1149)(6090,1088)
	(6081,1028)(6071,970)(6060,913)
	(6050,856)(6039,801)(6028,747)
	(6016,693)(6005,641)(5993,590)
	(5981,539)(5970,490)(5958,443)
	(5947,397)(5936,353)(5925,311)
	(5915,272)(5905,235)(5896,201)
	(5888,170)(5880,142)(5874,118)
	(5868,97)(5863,79)(5859,65)
	(5856,54)(5853,45)(5852,39)
	(5851,36)(5850,34)(5850,33)
\path(5850,2958)(5851,2959)(5852,2961)
	(5855,2964)(5859,2970)(5865,2977)
	(5872,2987)(5882,3000)(5892,3016)
	(5905,3033)(5919,3054)(5933,3076)
	(5949,3100)(5965,3127)(5981,3155)
	(5998,3185)(6014,3216)(6030,3249)
	(6045,3284)(6060,3321)(6075,3361)
	(6088,3402)(6101,3447)(6112,3494)
	(6123,3546)(6132,3601)(6139,3660)
	(6145,3723)(6149,3789)(6150,3858)
	(6149,3920)(6146,3982)(6142,4043)
	(6136,4103)(6128,4160)(6120,4216)
	(6110,4270)(6100,4322)(6089,4372)
	(6077,4421)(6065,4469)(6052,4515)
	(6038,4560)(6024,4604)(6010,4647)
	(5996,4689)(5981,4730)(5967,4769)
	(5953,4807)(5939,4843)(5925,4877)
	(5913,4908)(5901,4937)(5890,4963)
	(5881,4986)(5873,5005)(5866,5021)
	(5860,5034)(5856,5044)(5853,5050)
	(5851,5055)(5850,5057)(5850,5058)
\put(2700,3071){\ellipse{604}{604}}
\put(5550,2996){\ellipse{604}{604}}
\put(3600,333){\makebox(0,0)[lb]{\smash{{{\SetFigFont{5}{6.0}{\rmdefault}{\mddefault}{\updefault}$P_5$}}}}}
\put(6300,333){\makebox(0,0)[lb]{\smash{{{\SetFigFont{5}{6.0}{\rmdefault}{\mddefault}{\updefault}$P_3$}}}}}
\put(0,4233){\makebox(0,0)[lb]{\smash{{{\SetFigFont{8}{9.6}{\rmdefault}{\mddefault}{\updefault}$C_0^{-}$}}}}}
\end{picture}
}
 & $\cup$  &
\setlength{\unitlength}{0.0001in}
\begingroup\makeatletter\ifx\SetFigFont\undefined%
\gdef\SetFigFont#1#2#3#4#5{%
  \reset@font\fontsize{#1}{#2pt}%
  \fontfamily{#3}\fontseries{#4}\fontshape{#5}%
  \selectfont}%
\fi\endgroup%
{\renewcommand{\dashlinestretch}{30}
\begin{picture}(10908,5846)(0,-10)
\thicklines
\put(2325,1298){\blacken\ellipse{150}{150}}
\put(2325,1298){\ellipse{150}{150}}
\put(5475,698){\blacken\ellipse{150}{150}}
\put(5475,698){\ellipse{150}{150}}
\put(8325,398){\ellipse{600}{600}}
\put(8325,98){\blacken\ellipse{150}{150}}
\put(8325,98){\ellipse{150}{150}}
\path(75,698)(76,699)(78,700)
	(81,702)(85,705)(91,709)
	(99,714)(108,720)(119,728)
	(133,737)(148,747)(167,758)
	(187,771)(210,785)(235,800)
	(262,816)(292,833)(324,852)
	(358,871)(395,890)(433,911)
	(474,931)(517,952)(561,974)
	(608,995)(656,1017)(707,1038)
	(759,1059)(813,1080)(869,1101)
	(928,1121)(988,1140)(1051,1159)
	(1116,1178)(1184,1195)(1255,1212)
	(1328,1228)(1405,1243)(1485,1257)
	(1569,1270)(1657,1282)(1749,1292)
	(1845,1301)(1946,1309)(2052,1316)
	(2163,1320)(2279,1323)(2400,1324)
	(2526,1323)(2657,1321)(2792,1315)
	(2932,1308)(3075,1298)(3187,1289)
	(3300,1278)(3412,1266)(3524,1253)
	(3634,1239)(3742,1225)(3848,1209)
	(3951,1193)(4050,1176)(4146,1159)
	(4239,1141)(4328,1123)(4412,1105)
	(4493,1086)(4570,1068)(4643,1049)
	(4712,1030)(4777,1012)(4838,993)
	(4896,974)(4951,955)(5002,936)
	(5051,918)(5096,899)(5139,881)
	(5180,862)(5218,844)(5254,825)
	(5289,807)(5322,789)(5354,770)
	(5385,752)(5416,734)(5445,716)
	(5475,698)(5505,680)(5534,662)
	(5565,644)(5596,626)(5628,607)
	(5661,589)(5696,571)(5732,552)
	(5770,534)(5811,515)(5854,497)
	(5899,478)(5948,460)(5999,441)
	(6054,422)(6112,403)(6173,384)
	(6238,366)(6307,347)(6380,328)
	(6457,310)(6538,291)(6622,273)
	(6711,255)(6804,237)(6900,220)
	(6999,203)(7102,187)(7208,171)
	(7316,157)(7426,143)(7538,130)
	(7650,118)(7763,107)(7875,98)
	(8018,88)(8158,81)(8293,75)
	(8424,73)(8550,72)(8671,73)
	(8787,76)(8898,80)(9004,87)
	(9105,95)(9201,104)(9293,114)
	(9381,126)(9465,139)(9545,153)
	(9622,168)(9695,184)(9766,201)
	(9834,218)(9899,237)(9962,256)
	(10022,275)(10081,295)(10137,316)
	(10191,337)(10243,358)(10294,379)
	(10342,401)(10389,422)(10433,444)
	(10476,465)(10517,485)(10555,506)
	(10592,525)(10626,544)(10658,563)
	(10688,580)(10715,596)(10740,611)
	(10763,625)(10783,638)(10802,649)
	(10817,659)(10831,668)(10842,676)
	(10851,682)(10859,687)(10865,691)
	(10869,694)(10872,696)(10874,697)(10875,698)
\put(8500,2823){\makebox(0,0)[lb]{\smash{{{\SetFigFont{8}{9.6}{\rmdefault}{\mddefault}{\updefault}$C^+$}}}}}
\put(2358,700){\makebox(0,0)[lb]{\smash{{{\SetFigFont{5}{6.0}{\rmdefault}{\mddefault}{\updefault}$P_5$}}}}}
\put(5508,-150){\makebox(0,0)[lb]{\smash{{{\SetFigFont{5}{6.0}{\rmdefault}{\mddefault}{\updefault}$P_3$}}}}}
\end{picture}
}
\end{tabular}

The stable curve $C$ is represented by a point in some connected component $D$, which lies inside the stratum $M_3^{(5)}$. 
It decomposes as $C = C^+ \cup  C_0^- $, as pictured above, where  the subcurve $C^+$ is a $2$-pointed stable curve of genus $g^+=1$. 
\end{ex}

In general, if a pair $(C_0,P_1)$ is fixed, then there are numbers $g^+\le g$ and $m\le 4$, such that for any curve $C$ represented by a point of $D$ there exists a curve $C^+$, represented by a point of the moduli space $\overline{M}_{g^+,m}$ of pointed stable curves, with
\[ C \: = C^+ \cup C^-_0 .\]
Note  that while there is an  ordering of the marked points of $C_0^-$, which is induced by the embedding of $C_0^-$ into $C_0$, there is no such  ordering of the marked points of $C^+$. The process of glueing with $C_0^-$ determines the labels of the marked points of $C^+$ only up to certain permutations, which are induced by automorphisms of $C_0^-$. Thus, while the curve $C^+$ is not given as an $m$-pointed stable curve, it is uniquely determined as an $\mgamma$-pointed stable curve as introduced in \cite{Z2}, for a suitable group $\Gamma$. We will discuss this in detail in section \ref{sectauto}. 
 
In this way we can identify the subscheme $D$ with a known moduli space $M$, see lemma \ref{decompose}. Analogously, for the induced stratification of the moduli stack $\cmgbar$ of Deligne-Mumford stable curves, there is a relationship between the connected components ${\cal D}$ of the one-dimensional stratum of stable curves with $3g-4$ nodes, and corresponding moduli stacks ${\cal M}$. This relationship is our main topic, and the content of section \ref{sectmainthm}.

As we will see, the stacks ${\cal D}$ and ${\cal M}$ are in general not isomorphic to each other, even though their associated moduli spaces are. The reason for this is that the group of automorphisms of a curve $C$ is usually larger than that of a subcurve $C^+$. This difference does not affect the moduli spaces, and therefore does not affect $D$ and $M$, but it has an inpact on the associated moduli stacks. In fact, keeping track of automorphisms is what stacks are for. To recover an isomorphism, one essentially needs to enrich the stack ${\cal M}$ in order to take care of such  automorphisms, which act trivially on the curves $C^{\, \prime}$ represented by points of $M$, but which act nontrivially on its complement in $C^{\, \prime} \cup C_0^+$. 

For a general point of $D$ representing a curve $C$, its group of automorphisms $\Aut(C)$ splits in a natural way, see section \ref{sectauto} below. This allows us to modify the moduli stack $\cal M$ by taking the quotient with respect to one of the factors of $\Aut(C)$. Our main result is theorem \ref{MainThm}, stating that a dense open substack ${\cal D}^\times $ of ${\cal D}$ is isomorphic to a quotient of a dense open substack ${\cal M}^\times$ of ${\cal M}$. In particular, the stack ${\cal D}$ can generically be constructed very explicitely as a quotient stack, and incidentally as a quotient of another more accessible moduli stack. We obtain a surjective morphism 
\[ \calm^\times \: \longrightarrow \: \cald^\times ,\]
and show that is unramified and finite of degree equal to the order of a certain group of automorphisms.

If the group $\Aut(C_0)$ splits, too, then with certain additional technical effort it is possible to construct an analogous isomorphism, which extends over the boundary point representing the curve $C_0^+$, i.e. above the point in the closure of $D$ which  represents $C_0$. This is done in much detail in \cite{Z1}.  

For a special curve $C$, the group of automorphisms $\Aut(C)$ may not split. For this reason it is in general not possible to extend the above morphism to all of $\cal M$. As explained in \cite{Z1}, it is in particular impossible to generalize theorem \ref{MainThm} to a global quotient description of $\cald$, if there are points in $D$ representing curves with non-splitting automorphism groups. 

Our strategy to study the stacks ${\cal D}$ and ${\cal M}$ is to describe both as quotient stacks. We make use of the results of \cite{Z2}, where we generalized the construction of $\cmgbar$ by Edidin \cite{Ed} to moduli stacks of permutation classes of pointed stable curves. Using this quotient description, the relation between the stacks ${\cal M}$ and ${\cal D}$ is then given by a morphism between subschemes of suitable Hilbert schemes. It is straightforward to define a morphism $M\rightarrow D$ on the level of moduli spaces, using a process which Knudsen \cite{Kn} called clutching. However, its lifting to the level of Hilbert schemes requires some extra work, which is done in the sections \ref{sectglue} and \ref{sectdecompose}. In particular, at this stage the automorphism groups of stable curves play an important part. 

To illustrate our results, we conclude this paper with a complete analysis of the one-dimensional stratum in the boundary of the moduli stack $\cmgbar$ in the case of genus $g=3$. We compare our description to the classification Faber \cite{Fa} gave in the case of the moduli scheme $\overline{M}_3$.

\section{The moduli problem}\label{sectmodprob}

\begin{abschnitt}\em
The curves we are interested in are Deligne-Mumford stable curves of some genus $g \ge 3$ in the sense of Knudsen \cite{Kn}. He also introduced the more general notion of prestable curves. Here we will only be concerned with such prestable curves $f: C \rightarrow \Spec(k)$, where all connected components of $C$ are pointed stable curves themselves.

Throughout let one stable curve $f_0: C_0 \rightarrow \Spec(k)$ of genus $g\ge 3$ be fixed, which attains the maximal number $3g-3$ of nodes. We also choose an enumeration $P_1,\ldots,P_{3g-3}$ of the nodes. 
\end{abschnitt}

\begin{notation}\em\label{notc+}
Let $C_0^{+}$ denote the subscheme of the fixed curve $C_0$ which is the reduced union of all
irreducible components of $C_0$ containing the node $P_1$, together
with marked points $Q_1,\ldots,Q_m \in \{ P_2,\ldots,P_{3g-3}\}$, such
that $Q_i$ is not a singularity on $C_0^{+}$. Let $C_0^{-}$ denote the closure of the
complement of  $C_0^{+}$ in $C_0$, together with the same marked points $Q_1,\ldots,Q_m$. The set of all marked points shall be denoted by ${\cal Q}:=\{Q_1,\ldots,Q_m\}$. 
\end{notation}

\begin{abschnitt}\em
Note that the chosen enumeration of the nodes of $C_0$ distinguishes an ordering of the marked points of the subcurves.  
Thus the subcurve $f_0^{+}: C_0^{+} \rightarrow \Spec(k)$ is in a natural way an $m$-pointed stable curve of some genus
$g^{+}$. The subcurves $C_0^{+}$ and $C_0^{-}$
meet  transversally, and exactly in the points $Q_1,\ldots,Q_m$.  The curve $C_0^{-}$
may have more than one connected  
component, so it is an $m$-pointed prestable curve of some genus
$g^{-}$.  There is an equality  $g^{+}+g^{-}+m-1 =
g$. The
connected components of $C_0^{-}$ are again pointed stable curves, each with the maximal number of nodes possible.  
\end{abschnitt}

\begin{rk}\em\label{37}
Only very few types of $m$-pointed stable curves $C_0^{+}$ occur in
this way, and here is  a
complete list. Note that each irreducible component of $C_0$ is a rational curve with at most three nodes of $C_0$ lying on it, so we must have $m \le 4$. 

For $m=4$, there exist exactly three types of non-isomorphic nodal $4$-pointed
stable curves of genus $g^{+}=0$, each consisting of two rational
curves meeting each other transversally in one point. 

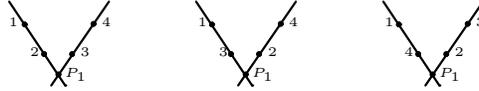
\begin{center}
\setlength{\unitlength}{0.00012in}
\begingroup\makeatletter\ifx\SetFigFont\undefined%
\gdef\SetFigFont#1#2#3#4#5{%
  \reset@font\fontsize{#1}{#2pt}%
  \fontfamily{#3}\fontseries{#4}\fontshape{#5}%
  \selectfont}%
\fi\endgroup%
{\renewcommand{\dashlinestretch}{30}
\begin{picture}(4341,3681)(0,-10)
\thicklines
\put(2155,505){\blacken\ellipse{150}{150}}
\put(2155,505){\ellipse{150}{150}}
\put(1533,1383){\blacken\ellipse{150}{150}}
\put(1533,1383){\ellipse{150}{150}}
\put(2733,1383){\blacken\ellipse{150}{150}}
\put(2733,1383){\ellipse{150}{150}}
\put(3655,2680){\blacken\ellipse{150}{150}}
\put(3655,2680){\ellipse{150}{150}}
\put(708,2658){\blacken\ellipse{150}{150}}
\put(708,2658){\ellipse{150}{150}}
\drawline(2433,33)(2433,33)
\path(33,3633)(2433,33)
\path(1833,33)(4308,3633)
\put(2433,408){\makebox(0,0)[lb]{\smash{{{\SetFigFont{5}{6.0}{\rmdefault}{\mddefault}{\updefault}$P_1$}}}}}
\put(33,2583){\makebox(0,0)[lb]{\smash{{{\SetFigFont{5}{6.0}{\rmdefault}{\mddefault}{\updefault}1}}}}}
\put(933,1308){\makebox(0,0)[lb]{\smash{{{\SetFigFont{5}{6.0}{\rmdefault}{\mddefault}{\updefault}2}}}}}
\put(4008,2583){\makebox(0,0)[lb]{\smash{{{\SetFigFont{5}{6.0}{\rmdefault}{\mddefault}{\updefault}4}}}}}
\put(3108,1308){\makebox(0,0)[lb]{\smash{{{\SetFigFont{5}{6.0}{\rmdefault}{\mddefault}{\updefault}3}}}}}
\put(-800,3308){\makebox(0,0)[lb]{\smash{{{\SetFigFont{5}{6.0}{\rmdefault}{\mddefault}{\updefault}{}}}}}}
\end{picture}
}
\hspace{5mm}
\setlength{\unitlength}{0.00012in}
\begingroup\makeatletter\ifx\SetFigFont\undefined%
\gdef\SetFigFont#1#2#3#4#5{%
  \reset@font\fontsize{#1}{#2pt}%
  \fontfamily{#3}\fontseries{#4}\fontshape{#5}%
  \selectfont}%
\fi\endgroup%
{\renewcommand{\dashlinestretch}{30}
\begin{picture}(4341,3681)(0,-10)
\thicklines
\put(2155,505){\blacken\ellipse{150}{150}}
\put(2155,505){\ellipse{150}{150}}
\put(1533,1383){\blacken\ellipse{150}{150}}
\put(1533,1383){\ellipse{150}{150}}
\put(2733,1383){\blacken\ellipse{150}{150}}
\put(2733,1383){\ellipse{150}{150}}
\put(3655,2680){\blacken\ellipse{150}{150}}
\put(3655,2680){\ellipse{150}{150}}
\put(708,2658){\blacken\ellipse{150}{150}}
\put(708,2658){\ellipse{150}{150}}
\drawline(2433,33)(2433,33)
\path(33,3633)(2433,33)
\path(1833,33)(4308,3633)
\put(2433,408){\makebox(0,0)[lb]{\smash{{{\SetFigFont{5}{6.0}{\rmdefault}{\mddefault}{\updefault}$P_1$}}}}}
\put(108,2583){\makebox(0,0)[lb]{\smash{{{\SetFigFont{5}{6.0}{\rmdefault}{\mddefault}{\updefault}1}}}}}
\put(1008,1308){\makebox(0,0)[lb]{\smash{{{\SetFigFont{5}{6.0}{\rmdefault}{\mddefault}{\updefault}3}}}}}
\put(3108,1308){\makebox(0,0)[lb]{\smash{{{\SetFigFont{5}{6.0}{\rmdefault}{\mddefault}{\updefault}2}}}}}
\put(4008,2583){\makebox(0,0)[lb]{\smash{{{\SetFigFont{5}{6.0}{\rmdefault}{\mddefault}{\updefault}4}}}}}
\put(-800,3308){\makebox(0,0)[lb]{\smash{{{\SetFigFont{5}{6.0}{\rmdefault}{\mddefault}{\updefault}{}}}}}}
\end{picture}
}
\hspace{5mm}
\setlength{\unitlength}{0.00012in}
\begingroup\makeatletter\ifx\SetFigFont\undefined%
\gdef\SetFigFont#1#2#3#4#5{%
  \reset@font\fontsize{#1}{#2pt}%
  \fontfamily{#3}\fontseries{#4}\fontshape{#5}%
  \selectfont}%
\fi\endgroup%
{\renewcommand{\dashlinestretch}{30}
\begin{picture}(4341,3681)(0,-10)
\thicklines
\put(2155,505){\blacken\ellipse{150}{150}}
\put(2155,505){\ellipse{150}{150}}
\put(1533,1383){\blacken\ellipse{150}{150}}
\put(1533,1383){\ellipse{150}{150}}
\put(2733,1383){\blacken\ellipse{150}{150}}
\put(2733,1383){\ellipse{150}{150}}
\put(3655,2680){\blacken\ellipse{150}{150}}
\put(3655,2680){\ellipse{150}{150}}
\put(708,2658){\blacken\ellipse{150}{150}}
\put(708,2658){\ellipse{150}{150}}
\drawline(2433,33)(2433,33)
\path(33,3633)(2433,33)
\path(1833,33)(4308,3633)
\put(2433,408){\makebox(0,0)[lb]{\smash{{{\SetFigFont{5}{6.0}{\rmdefault}{\mddefault}{\updefault}$P_1$}}}}}
\put(108,2583){\makebox(0,0)[lb]{\smash{{{\SetFigFont{5}{6.0}{\rmdefault}{\mddefault}{\updefault}1}}}}}
\put(933,1308){\makebox(0,0)[lb]{\smash{{{\SetFigFont{5}{6.0}{\rmdefault}{\mddefault}{\updefault}4}}}}}
\put(3108,1308){\makebox(0,0)[lb]{\smash{{{\SetFigFont{5}{6.0}{\rmdefault}{\mddefault}{\updefault}2}}}}}
\put(4008,2583){\makebox(0,0)[lb]{\smash{{{\SetFigFont{5}{6.0}{\rmdefault}{\mddefault}{\updefault}3}}}}}
\put(-800,3308){\makebox(0,0)[lb]{\smash{{{\SetFigFont{5}{6.0}{\rmdefault}{\mddefault}{\updefault}{}}}}}}
\end{picture}
}
\end{center}

The curves pictured above are represented by the three points which lie in the boundary of the
moduli space $\msbar_{0,4}$. If the condition on the existence of a
node at $P_1$ is dropped, we 
obtain the full moduli space $\msbar_{0,4}$.

If $m=2$, then  there exist two types of non-isomorphic nodal $2$-pointed
stable curves of genus $g^{+}=1$, both with two nodes and two rational irreducible components, as pictured below.

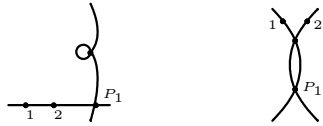
\begin{center}
\setlength{\unitlength}{0.00012in}
\begingroup\makeatletter\ifx\SetFigFont\undefined%
\gdef\SetFigFont#1#2#3#4#5{%
  \reset@font\fontsize{#1}{#2pt}%
  \fontfamily{#3}\fontseries{#4}\fontshape{#5}%
  \selectfont}%
\fi\endgroup%
{\renewcommand{\dashlinestretch}{30}
\begin{picture}(4416,5106)(0,-10)
\thicklines
\put(3258,2996){\ellipse{604}{604}}
\path(3558,3033)(3558,3032)(3559,3031)
	(3561,3028)(3564,3023)(3568,3017)
	(3574,3008)(3581,2996)(3589,2983)
	(3598,2967)(3609,2948)(3621,2927)
	(3633,2903)(3647,2878)(3661,2850)
	(3675,2821)(3689,2789)(3704,2756)
	(3718,2721)(3732,2684)(3746,2644)
	(3760,2603)(3773,2559)(3785,2513)
	(3797,2463)(3808,2410)(3818,2353)
	(3828,2293)(3836,2227)(3844,2157)
	(3850,2083)(3854,2004)(3857,1920)
	(3858,1833)(3857,1761)(3855,1689)
	(3852,1617)(3848,1545)(3843,1475)
	(3837,1407)(3831,1340)(3823,1275)
	(3815,1211)(3807,1149)(3798,1088)
	(3789,1028)(3779,970)(3768,913)
	(3758,856)(3747,801)(3736,747)
	(3724,693)(3713,641)(3701,590)
	(3689,539)(3678,490)(3666,443)
	(3655,397)(3644,353)(3633,311)
	(3623,272)(3613,235)(3604,201)
	(3596,170)(3588,142)(3582,118)
	(3576,97)(3571,79)(3567,65)
	(3564,54)(3561,45)(3560,39)
	(3559,36)(3558,34)(3558,33)
\path(3558,2958)(3559,2959)(3560,2961)
	(3563,2964)(3567,2970)(3573,2977)
	(3580,2987)(3590,3000)(3600,3016)
	(3613,3033)(3627,3054)(3641,3076)
	(3657,3100)(3673,3127)(3689,3155)
	(3706,3185)(3722,3216)(3738,3249)
	(3753,3284)(3768,3321)(3783,3361)
	(3796,3402)(3809,3447)(3820,3494)
	(3831,3546)(3840,3601)(3847,3660)
	(3853,3723)(3857,3789)(3858,3858)
	(3857,3920)(3854,3982)(3850,4043)
	(3844,4103)(3836,4160)(3828,4216)
	(3818,4270)(3808,4322)(3797,4372)
	(3785,4421)(3773,4469)(3760,4515)
	(3746,4560)(3732,4604)(3718,4647)
	(3704,4689)(3689,4730)(3675,4769)
	(3661,4807)(3647,4843)(3633,4877)
	(3621,4908)(3609,4937)(3598,4963)
	(3589,4986)(3581,5005)(3574,5021)
	(3568,5034)(3564,5044)(3561,5050)
	(3559,5055)(3558,5057)(3558,5058)
\put(3558,2958){\blacken\ellipse{150}{150}}
\put(3558,2958){\ellipse{150}{150}}
\put(3783,708){\blacken\ellipse{150}{150}}
\put(3783,708){\ellipse{150}{150}}
\put(783,708){\blacken\ellipse{150}{150}}
\put(783,708){\ellipse{150}{150}}
\put(1983,708){\blacken\ellipse{150}{150}}
\put(1983,708){\ellipse{150}{150}}
\path(33,708)(4383,708)
\put(4083,1008){\makebox(0,0)[lb]{\smash{{{\SetFigFont{5}{6.0}{\rmdefault}{\mddefault}{\updefault}$P_1$}}}}}
\put(783,108){\makebox(0,0)[lb]{\smash{{{\SetFigFont{5}{6.0}{\rmdefault}{\mddefault}{\updefault}$1$}}}}}
\put(1983,108){\makebox(0,0)[lb]{\smash{{{\SetFigFont{5}{6.0}{\rmdefault}{\mddefault}{\updefault}$2$}}}}}
\end{picture}
}
\hspace*{1.5cm}
\setlength{\unitlength}{0.00012in}
\begingroup\makeatletter\ifx\SetFigFont\undefined%
\gdef\SetFigFont#1#2#3#4#5{%
  \reset@font\fontsize{#1}{#2pt}%
  \fontfamily{#3}\fontseries{#4}\fontshape{#5}%
  \selectfont}%
\fi\endgroup%
{\renewcommand{\dashlinestretch}{30}
\begin{picture}(2016,4881)(0,-10)
\thicklines
\put(1008,3483){\blacken\ellipse{150}{150}}
\put(1008,3483){\ellipse{150}{150}}
\put(1008,1383){\blacken\ellipse{150}{150}}
\put(1008,1383){\ellipse{150}{150}}
\put(1533,4308){\blacken\ellipse{150}{150}}
\put(1533,4308){\ellipse{150}{150}}
\put(483,4308){\blacken\ellipse{150}{150}}
\put(483,4308){\ellipse{150}{150}}
\path(33,4833)(34,4832)(35,4831)
	(37,4829)(41,4824)(47,4819)
	(55,4811)(65,4800)(77,4787)
	(92,4772)(109,4754)(129,4733)
	(151,4709)(176,4682)(203,4653)
	(232,4621)(263,4586)(296,4549)
	(330,4510)(366,4468)(403,4425)
	(441,4379)(480,4332)(519,4283)
	(558,4233)(597,4181)(637,4128)
	(676,4074)(714,4018)(753,3961)
	(790,3903)(827,3844)(863,3783)
	(898,3720)(932,3656)(964,3590)
	(996,3523)(1026,3453)(1055,3381)
	(1082,3307)(1108,3230)(1131,3151)
	(1153,3070)(1173,2985)(1190,2899)
	(1205,2809)(1217,2717)(1226,2624)
	(1231,2529)(1233,2433)(1231,2337)
	(1226,2242)(1217,2149)(1205,2057)
	(1190,1967)(1173,1881)(1153,1796)
	(1131,1715)(1108,1636)(1082,1559)
	(1055,1485)(1026,1413)(996,1343)
	(964,1276)(932,1210)(898,1146)
	(863,1083)(827,1022)(790,963)
	(753,905)(714,848)(676,792)
	(637,738)(597,685)(558,633)
	(519,583)(480,534)(441,487)
	(403,441)(366,398)(330,356)
	(296,317)(263,280)(232,245)
	(203,213)(176,184)(151,157)
	(129,133)(109,112)(92,94)
	(77,79)(65,66)(55,55)
	(47,47)(41,42)(37,37)
	(35,35)(34,34)(33,33)
\path(1983,4833)(1982,4832)(1981,4831)
	(1979,4829)(1975,4824)(1969,4819)
	(1961,4811)(1951,4800)(1939,4787)
	(1924,4772)(1907,4754)(1887,4733)
	(1865,4709)(1840,4682)(1813,4653)
	(1784,4621)(1753,4586)(1720,4549)
	(1686,4510)(1650,4468)(1613,4425)
	(1575,4379)(1536,4332)(1497,4283)
	(1458,4233)(1419,4181)(1379,4128)
	(1340,4074)(1302,4018)(1263,3961)
	(1226,3903)(1189,3844)(1153,3783)
	(1118,3720)(1084,3656)(1052,3590)
	(1020,3523)(990,3453)(961,3381)
	(934,3307)(908,3230)(885,3151)
	(863,3070)(843,2985)(826,2899)
	(811,2809)(799,2717)(790,2624)
	(785,2529)(783,2433)(785,2337)
	(790,2242)(799,2149)(811,2057)
	(826,1967)(843,1881)(863,1796)
	(885,1715)(908,1636)(934,1559)
	(961,1485)(990,1413)(1020,1343)
	(1052,1276)(1084,1210)(1118,1146)
	(1153,1083)(1189,1022)(1226,963)
	(1263,905)(1302,848)(1340,792)
	(1379,738)(1419,685)(1458,633)
	(1497,583)(1536,534)(1575,487)
	(1613,441)(1650,398)(1686,356)
	(1720,317)(1753,280)(1784,245)
	(1813,213)(1840,184)(1865,157)
	(1887,133)(1907,112)(1924,94)
	(1939,79)(1951,66)(1961,55)
	(1969,47)(1975,42)(1979,37)
	(1981,35)(1982,34)(1983,33)
\put(1308,1308){\makebox(0,0)[lb]{\smash{{{\SetFigFont{5}{6.0}{\rmdefault}{\mddefault}{\updefault}$P_1$}}}}}
\put(-150,4000){\makebox(0,0)[lb]{\smash{{{\SetFigFont{5}{6.0}{\rmdefault}{\mddefault}{\updefault}$1$}}}}}
\put(1800,4000){\makebox(0,0)[lb]{\smash{{{\SetFigFont{5}{6.0}{\rmdefault}{\mddefault}{\updefault}$2$}}}}}
\end{picture}
}
\end{center}

Both of them are represented by isolated points in the two-dimensional
moduli space $\msbar_{1,2}$. The locus of curves in $\msbar_{1,2}$,
which preserve the node which is not $P_1$, is a one-dimensional subscheme $\Delta_0$. In fact, $\Delta_0$ is given by that connected component of the boundary stratum $M_{1,2}^{(1)}$ of $\msbar_{1,2}$, which represents irreducible singular stable curves.  

The only instance with $m=1$ is a curve of arithmetic genus $g^+=1$ with
one node at $P_1$, i.e. the one point in the boundary of $\msbar_{1,1}$.

\begin{center}
\setlength{\unitlength}{0.00012in}
\begingroup\makeatletter\ifx\SetFigFont\undefined%
\gdef\SetFigFont#1#2#3#4#5{%
  \reset@font\fontsize{#1}{#2pt}%
  \fontfamily{#3}\fontseries{#4}\fontshape{#5}%
  \selectfont}%
\fi\endgroup%
{\renewcommand{\dashlinestretch}{30}
\begin{picture}(5091,1314)(0,-10)
\thicklines
\put(2996,975){\ellipse{604}{604}}
\path(3033,675)(3032,675)(3031,674)
	(3028,672)(3023,669)(3017,665)
	(3008,659)(2996,652)(2983,644)
	(2967,635)(2948,624)(2927,612)
	(2903,600)(2878,586)(2850,572)
	(2821,558)(2789,544)(2756,529)
	(2721,515)(2684,501)(2644,487)
	(2603,473)(2559,460)(2513,448)
	(2463,436)(2410,425)(2353,415)
	(2293,405)(2227,397)(2157,389)
	(2083,383)(2004,379)(1920,376)
	(1833,375)(1761,376)(1689,378)
	(1617,381)(1545,385)(1475,390)
	(1407,396)(1340,402)(1275,410)
	(1211,418)(1149,426)(1088,435)
	(1028,444)(970,454)(913,465)
	(856,475)(801,486)(747,497)
	(693,509)(641,520)(590,532)
	(539,544)(490,555)(443,567)
	(397,578)(353,589)(311,600)
	(272,610)(235,620)(201,629)
	(170,637)(142,645)(118,651)
	(97,657)(79,662)(65,666)
	(54,669)(45,672)(39,673)
	(36,674)(34,675)(33,675)
\path(2958,675)(2959,674)(2961,673)
	(2964,670)(2970,666)(2977,660)
	(2987,653)(3000,643)(3016,633)
	(3033,620)(3054,606)(3076,592)
	(3100,576)(3127,560)(3155,544)
	(3185,527)(3216,511)(3249,495)
	(3284,480)(3321,465)(3361,450)
	(3402,437)(3447,424)(3494,413)
	(3546,402)(3601,393)(3660,386)
	(3723,380)(3789,376)(3858,375)
	(3920,376)(3982,379)(4043,383)
	(4103,389)(4160,397)(4216,405)
	(4270,415)(4322,425)(4372,436)
	(4421,448)(4469,460)(4515,473)
	(4560,487)(4604,501)(4647,515)
	(4689,529)(4730,544)(4769,558)
	(4807,572)(4843,586)(4877,600)
	(4908,612)(4937,624)(4963,635)
	(4986,644)(5005,652)(5021,659)
	(5034,665)(5044,669)(5050,672)
	(5055,674)(5057,675)(5058,675)
\put(2958,675){\blacken\ellipse{150}{150}}
\put(2958,675){\ellipse{150}{150}}
\put(1008,450){\blacken\ellipse{150}{150}}
\put(1008,450){\ellipse{150}{150}}
\put(2808,-200){\makebox(0,0)[lb]{\smash{{{\SetFigFont{5}{6.0}{\rmdefault}{\mddefault}{\updefault}$P_1$}}}}}
\put(1008,-200){\makebox(0,0)[lb]{\smash{{{\SetFigFont{5}{6.0}{\rmdefault}{\mddefault}{\updefault}$1$}}}}}
\end{picture}
}
\end{center}

If the condition on the existence of a node at $P_1$ is dropped, we
obtain the full moduli space $\msbar_{1,1}$. 
\end{rk}

\begin{notation}\em\label{mprimed}
For the integers $m=1,2$ or $4$ which occur as numbers of marked points of $C_0^-$, we define
\[ M_{g^{+},m}' \: := \: \left\{ \begin{array}{ll}
M_{1,1} & \Text{for } m = 1;\\
\Delta_0& \Text{for } m=2;\\
M_{0,4} & \Text{for } m = 4.
\end{array} \right. \]
Here, as usual $M_{g^{+},m}$ denotes the open subscheme of $\msbar_{g^{+},m}$, parametrizing nonsingular stable curves. By
$\Delta_0 \subset \msbar_{1,2}$ we denote the reduced subscheme given by the locus of irreducible singular stable curves with exactly one node.
Note  that  
in each case the closure of $M_{g^{+},m}'$ is isomorphic to $ \PP^1$.

The pair $(C_0,P_1)$, which we fixed above, induces naturally a morphism 
\[ \theta : \quad M_{g^{+},m}' \: \longrightarrow \: \msbar_{g} \]
as follows. For a point $[C] \in M_{g^{+},m}' $, represented by an $m$-pointed stable curve $f: C\rightarrow \Spec(k)$, its image $\theta([C])$ in $ \msbar_{g}$ is defined as the point representing the stable curve of genus $g$ obtained by glueing $C$ and $C_0^-$ in their respective marked points. 

Since the scheme $M_{g^{+},m}'$ is connected, this morphism surjects onto one connected component of the boundary stratum $M_{g}^{(3g-4)}$. By construction, the point representing $C_0$ is contained in the closure of the image of $\theta$. 
\end{notation}

\begin{defi}\em
Let a stable curve $f_0: C_0\rightarrow \Spec(k)$ with $3g-3$ nodes and a node $P_1$ of $C_0$ be given. We denote by
\[ \dpo \: \subset \: \msbar_g \]
the connected component of the boundary stratum $M_{g}^{(3g-4)}$, which is the reduced image of the above morphism $\theta$.  
\end{defi}

Decomposing a stable curve $f: C\rightarrow \Spec(k)$ by separating it at such nodes, where two irreducible components meet, is fairly straightforward. When we are dealing with families of curves, i.e. stable curves $f: C\rightarrow S$ over a general base scheme $S$, the situation becomes more subtle. To be able to ``cut'' along the nodes of $C$, the nodes of the fibres must be given by sections $\sigma: S \rightarrow C$, and in general this is not the case. 

\begin{lemma}\label{sectnod}
Let $f: C\rightarrow S$ be a stable curve of genus $g$, with reduced and irreducible base $S$. Suppose that the number of nodes of the fibre over a general closed point is $k$, for some $1\le k \le 3g-3$. Then there exists an \'etale covering $u: S'\rightarrow S$, where $S'$ is connected, and such that the pullback $f': C'\rightarrow S'$ of $f:C\rightarrow S$ to $S'$ admits $k$ disjoint sections of nodes.  
\end{lemma}

\proof
Note that the assumption on a general fibre implies that for an arbitrary closed point $s\in S$, the fibre $C_s$ has at least $k$ nodes. Consider the subscheme $N \subset C$, which is defined as the reduced locus of all points of $C$ which are nodes of fibres of $f:C\rightarrow S$. Formally, it can be constructed as the locus of those points, where the sheaf $\Omega_{C/S}^1$ is not free. Because of the flatness of $f:C\rightarrow S$, these points are precisely the nodes of the fibres. Define $R$ as the reduced union of all those irreducible components of $N$, which surject onto $S$.

Let $s\in S$ be a general point, and let $n\in C_s$ be one of the $k$  nodes over $s$. Since all neighbouring fibres have $k$ nodes as well, we must have $n\in R$. Therefore, the fibre of the restricted morphism $f|R: R\rightarrow S$ over $s$ consists of exactly $k$ necessarily simple points, which are precisely all the nodes over $s$. Clearly no fibre of $f|R$ may consist of more than $k$ points. Recall that closed points of $R$ are nodes of fibres of $f:C\rightarrow S$, so they are ordinary double points. There is no deformation, which splits one ordinary double point into two, so there can also be no fibre of $f|R$ consisting of less than $k$ points. In other words, each fibre of $f|R: R\rightarrow S$ consists of  exactly $k$ simple points, and  therefore the morphism $f|R:R\rightarrow S$ is flat, and an \'etale covering of $S$. Since $S$ is connected, it holds that for any connected component $S_1$ of $R$ the restriction $g:= f|S_1: S_1\rightarrow S$ is an \'etale covering, too. 

We now define $f_1:C_1\rightarrow S_1$ as the pullback of $f:C\rightarrow S$ to $S_1$. In particular, the stable curve $f_1:C_1\rightarrow S_1$ admits a tautological section $\sigma_1: S_1\rightarrow C_1$, such that for all closed points $s\in S_1$, the point $\sigma_1(s)$ is a node of the fibre $C_{1,s}$. 

Note that $\sigma(S_1)$ is a connected component of the pullback $g^\ast R\subset C_1$, since the restricted morphism $g^\ast R \rightarrow S_1$ is \'etale. Repeating the above construction with the complement $R_1 := g^\ast R \smallsetminus \sigma(S_1)$ as a subscheme of $C_1$, we obtain an \'etale morphisms $R_1\rightarrow S_1$ which is of degree $k-1$. Proceeding inductively, we define a sequence of \'etale coverings $S_{i+1} \rightarrow S_{i}$, for $i=1,\ldots,k-1$, such that the combined pullback $f': C' \rightarrow S'$ of $f:C\rightarrow S$ to $S':= S_{k+1}$ admits $k$ sections of nodes.
\ebew

\section{Automorphisms of decomposing stable curves}\label{sectauto}

The automorphisms of stable curves will play an important part in our considerations below. We need to relate the group of automorphisms of a stable curve to the automorphism groups of its subcurves. In particular we need to understand the group $\Aut(C_0^-)$.

\begin{notation}\em\label{not31}
Let $\CQ:= \{Q_1,\ldots,Q_m\}$ denote the set of marked points of an $m$-pointed curve $f: C\rightarrow \Spec(k)$, like for example $C_0^-$ or $C_0^+$. We define 
$ \Aut(C;\CQ) $ as 
the group of all those automorphisms of $C$ as a scheme over $\Spec(k)$, which map the set $\CQ$ to itself. As usual, $\Aut(C)$ denotes the automorphism group of $C$ as an $m$-pointed curve. 

Let $\Sigma_m$ denote the permutation group acting on the labels of the $m$ marked points of $C$ contained in $\CQ$.
For any $\phi \in  \Aut(C;\CQ)$ there is a unique permutation $\gamma \in \Sigma_m$ such that $\phi(Q_i) = Q_{\gamma(i)}$
holds for all $i=1,\ldots,m$. Thus there is a natural homomorphism of groups
$ \Aut(C;\CQ) \rightarrow \Sigma_m $. 
We define  $\Gamma(C;\CQ)$ as the image of this homomorphism. In particular there is an exact sequence
\[ \id \rightarrow \Aut(C)\rightarrow \Aut(C;\calq)  \rightarrow \Gamma(C;\CQ)\rightarrow \id .\]

Let $\Gamma \subset \Sigma_m$ be a subgroup. 
By definition, an $\mgamma$-pointed stable curve over $\Spec(k)$ is an equivalence class of $m$-pointed stable curves $f: C\rightarrow \Spec(k)$ with respect to the action of $\Gamma$ on the labels of the marked points, as introduced in  \cite{Z2}. It is also called a permutation class of stable curves. The group of automorphisms of the $\mgamma$-pointed stable curve represented by $C$ is denoted by $\Aut_\Gamma(C;\calq)$. It is the subgroup of those elements  $\phi\in\Aut(C;\calq)$, for which  there exists a permutation $\gamma \in \Gamma$ such that 
$ \phi(Q_i) = Q_{\gamma(i)} $
holds for all $i=1,\ldots,m$. 
\end{notation}

\begin{lemma}\label{newsplit}
Let $f^+: C^+\rightarrow \Spec(k)$ be a stable curve of genus $g^+$ with $m$ marked points, which is represented by a point $[C^+]\in M_{g^{+},m}'$. Let $f: C\rightarrow \Spec(k)$ be the stable curve of genus $g$, which is obtained by glueing $C^+$ with $C_0^-$ in the $m$ marked points. Then there is an exact sequence of groups
\[ \id \rightarrow \Aut(C_0^-) \rightarrow \Aut(C) \rightarrow \Aut_\Gamma(C^+;\calq) \rightarrow \id ,\]
where  $\Gamma= \Gamma(C_0^-;\calq)$. 
\end{lemma}

\proof
By definition, any automorphism in $\Aut(C_0^-)$ fixes each one of the marked points of $C_0^-$, and these marked points are the glueing points of $C_0^-$ with $C^+$. By trivial extension, the group $\Aut(C_0^-)$ becomes a subgroup of $\Aut(C)$. Conversely, since $[C^+]\in M_{g^{+},m}'$, the curve $C^+$  is the unique connected subcurve of $C$ of genus $g^+$, which has not the maximal number of nodes as an $m$-pointed stable curve. Hence any automorphism of $C$ restricts to an automorphism on $C_0^-$ and to an automorphism on $C^+$, which both are not necessarily fixing the marked points. However, in this way only such permutations of marked points occur on $C^+$, which are also effected by automorphisms on $C_0^-$. Thus by definition, an automorphism of $C$ restricted to $C^+$ is an automorphism of $C^+$ considered as an $\mgamma$-pointed stable curve, i.e. an element of $\Aut_\Gamma(C^+;\calq)$. 

All automorphisms in $\Aut_\Gamma(C^+;\calq)$ occur in this way. To see this, let $\rho^+ \in \Aut_\Gamma(C^+;\calq)$ be an automorphism of $C^+$ as an  $\mgamma$-pointed stable curve. It induces a permutation $\gamma$ of the marked points, with $\gamma\in \Gamma$. By the definition of $\Gamma$, there exists an automorphism $\rho^-\in\Aut(C_0^-;\calq)$ inducing $\gamma$. Glueing the automorphisms $\rho^+$ and $\rho^-$ defines an automorphism $\rho$ on $C$, which restricts to $\rho^+$. Thus we have shown that there exists an exact sequence of groups, as claimed.
\ebew

\begin{rk}\em\label{newsplitrk}
For a general point $[C^+]\in M_{g^{+},m}'$ as in lemma \ref{newsplit}, the exact sequence of the lemma  splits in a natural way. This can be seen as follows. If $m=1$, then clearly $\Aut_\Gamma(C^+;\calq) = \Aut(C^+)$. If $m = 2$ or $m=4$, then an examination of the list of possibilities for curves $C^+$ from  remark \ref{37} shows that automorphisms of $C^+$, which permute the marked points, exist only if the marked points are in a very special position. 
Thus, if $[C^+]\in {M}_{g^+,m}'$ is a general point, then again $\Aut_\Gamma(C^+;\calq) = \Aut(C^+)$. Now a splitting homomorphism $\Aut_\Gamma(C^+;\calq) \rightarrow \Aut(C)$ can be defined by trivial extension.
\end{rk}

\begin{lemma}\label{decompose}
Let $f : C\rightarrow \Spec(k)$ be a stable curve of genus $g\ge 3$ represented by a point $[C]\in \dpo$. Let $\Gamma = \Gamma(C_0^-;\calq)$. Then there exists a uniquely determined $\mgamma$-pointed stable curve $C^+$, with $m\in\{1,2,4\}$, such that $C$ is obtained from glueing a representative of $C^+$ and the $m$-pointed prestable curve $C_0^-$ in the corresponding marked points. 
In particular, there is an isomorphism of schemes 
\[ \dpo \cong  {M}_{g^+,m}'/\Gamma .\]
\end{lemma}

Let us remark that the closure of $\dpo$ in $\mgbar$ is always isomorphic to $\PP^1$, see \cite{Z1}. 

\proof
By the definition of $\dpo$, there exists a decomposition of the curve $C$ into an $m$-pointed stable subcurve $C^+$ and a subcurve isomorphic to $C_0^-$. The labels of the marked points of $C^+$ are determined only up to certain permutations, which are induced by automorphisms of $C_0^-$, i.e. $C^+$  is uniquely determined as an $\mgamma$-pointed stable curve. 
\ebew

\begin{defi}\em\label{def_split_d}
We define the subscheme
\[ \dpsplit \: \subset \:\dpo \]
as the nonempty reduced open subscheme, which is the locus of points $[C]$ para\-metrizing those curves  $f:C\rightarrow \Spec(k)$, whose groups of automorphisms split naturally as
$ \Aut(C) = \Aut_\Gamma(C^+;\calq)\times \Aut(C_0^-) $, 
where $f^+: C^+\rightarrow \Spec(k)$ denotes the unique subcurve as in lemma \ref{decompose}.
\end{defi} 

\section{The glueing morphism}\label{sectglue}

Our objective in this section is to construct a lifting of the natural morphism $\theta : {M}'_{g^+,m} \rightarrow \overline{M}_g$ from paragraph \ref{mprimed} generically to a morphism between subschemes of Hilbert schemes. 

\begin{notation}\em\label{not32}
Let $f: C\rightarrow S$ be an $m$-pointed  stable curve of genus $g$ over some base scheme $S$, such that  $2g-2+m > 0$. Let $S_1,\ldots ,S_{m}$ denote the divisors of marked points of $C$, and let $\omega_{C/S}$ be the canonical sheaf. If $n\in{\mathbb N}$ is chosen large enough, then 
 the sheaf $\left(\omega_{C/S}(S_1+\ldots+S_{m})\right)^{\otimes n}$ is
relatively very ample, and $f_\ast\left(\left(\omega_{C/S}(S_1+\ldots+S_{m})\right)^{\otimes n}\right)$ is locally free of
rank $(2g-2+m)n -g+1$, compare \cite{Kn}. 

Consider the Hilbert scheme $\Hilb_{{\PP^N}}^{{P_{g,n,m}}}$ of curves $ f: C \rightarrow \Spec(k)$ embedded in $\PP^N$ with Hilbert polynomial $P_{g,n,m}(t) := (2g-2+m)nt-g+1$, where  $N :=(2g-2+m)n -g$. We denote by
\[ \overline{H}_{g,n,m} \: \: \subset \: \: \Hilb_{{\PP^N}}^{{P_{g,n,m}}} \times (\PP^N)^m\]
the quasi-projective reduced subscheme, which is the locus of $n$-canonically embedded $m$-pointed stable curves.  

The standard construction of the moduli space $\overline{M}_{g,m}$ of $m$-pointed stable curves of genus $g$ is its construction as the GIT-quotient of $\overline{H}_{g,n,m}$ with respect to the natural $\PGL(N+1)$-action, compare for example \cite{DM} or \cite{Gs}.

If $\Gamma \subset \Sigma_m$ is a subgroup of the group of permutations on the labels of the $m$ marked points, then the moduli space of $\mgamma$-pointed stable curves is the quotient $\overline{M}_{g,\mgamma}:= \overline{M}_{g,m}/\Gamma$, see \cite{Z2}. Equivalently, it can be constructed as the GIT-quotient of $\overline{H}_{g,n,\mgamma}:= \overline{H}_{g,n,m}/\Gamma$ with respect to the natural $\PGL(N+1)$-action. 
\end{notation}

\begin{abschnitt}\em
As before, let $f_0:C_0\rightarrow \Spec(k)$ denote the fixed stable curve of genus $g\ge 3$ with $3g-3$ nodes, and let $P_1,\ldots,P_{3g-3}$ be an enumeration of its nodes. Let $C_0^-$ denote the unique $m$-pointed prestable curve of genus $g^-$, which is the reduced union of all of the irreducible components of $C_0$ which do not contain $P_1$. The closure of the complement of $C_0^-$ in $C_0$ is denoted by $C_0^+$. It can be viewed as an $m$-pointed stable curve of some genus $g^+$. 
\end{abschnitt}

\begin{notation}\em
Let $\Gamma := \Gamma(C_0^-;\CQ)$. 
For $n\in{\mathbb N}$ sufficiently large let $\overline{H}_{g,n,0}$ and $\overline{H}_{g^+,n,\mgamma}$ denote the subschemes of the Hilbert schemes as in \ref{not32}, with the moduli spaces $\overline{M}_{g}$ and $\overline{M}_{g^+,\mgamma}$ as their respective GIT-quotients. Let ${H}_{g^+,n,m}'$ denote the reduced preimage of ${M}_{g^+,m}'$ in $\overline{H}_{g^+,n,m}$, and we define ${H}_{g^+,n,\mgamma}' := {H}_{g^+,n,m}'/\Gamma$, as well as ${M}_{g^+,\mgamma}' := {M}_{g^+,m}'/\Gamma$.
\end{notation}

\begin{abschnitt}\em\label{const_part_1}
{\em Construction step 1.} 
Consider the universal embedded $m$-pointed stable curve $u^{+}: {\cal C}_{{g^{+},n,m}} \rightarrow
\hilb{g^{+},n,m}$, which is the restriction of the universal curve on the full Hilbert scheme. We glue this along the $m$ given sections to the
trivial $m$-pointed prestable curve $C_0^{-}\times \hilb{g^{+},n,m}
\rightarrow \hilb{g^{+},n,m}$. This produces a stable curve $u_0: {\cal
C}_0\rightarrow  \hilb{g^{+},n,m}$ of genus
$g$, compare \cite{Kn}. Knudsen calls this procedure of glueing along a pair of disjoint sections ``clutching''. 

Let  $i : {\cal C}_{{g^{+},n,m}} \hookrightarrow {\cal C}_0 $ denote the inclusion morphism. 
Consider the restriction of the canonical invertible dualizing sheaf $\omega_{{\cal C}_0/\hilb{g^{+},n,m}}$ to a fibre $C_h \rightarrow \Spec(k)$ of $u_0: {\cal
C}_0\rightarrow  \hilb{g^{+},n,m}$, for some closed point $h\in \hilb{g^{+},n,m}$. Let $n : \hat{C}_h \rightarrow C_h$ be the normalization of $C_h$, and for a node $P_i\in C_h$ let $A_i$ and $B_i$ denote its preimages in $\hat{C}_h$. 
The restriction of $\omega_{{\cal C}_0/\hilb{g^{+},n,m}}$ to $C_h$ can be described as the sheaf of those $1$-forms on $\hat{C}_h$ which have at most simple poles in $A_i$ and $B_i$, and such that the sum of the residues in $A_i$ and $B_i$ is zero, compare for example \cite{Kn}. 

By construction, the curve $u_0: {\cal
C}_0\rightarrow \hilb{g^{+},n,m}$ comes with $m$ disjoint sections of nodes, which are given by the sections of glueing. Let $n_0 : \hat{C}_h' \rightarrow C_h$ denote the partial normalization of $C_h$, in the sense that it is the normalization map near those  nodes resulting from the glueing, and an isomorphism elsewhere. In particular, $\hat{C}_h'$ preserves exactly those nodes, which are not glueing nodes.  Thus $\hat{C}_h'$ decomposes into two disjoint subcurves, one of them isomorphic to $C_0^-$, the other one isomorphic to $C^+_h$, where $C_h^+$ is the fibre of the universal curve $u^{+}: {\cal C}_{{g^{+},n,m}} \rightarrow
\hilb{g^{+},n,m}$ over the point $h$. This induces a morphism $\iota: C_h^+ \rightarrow \hat{C}_h'$, such that $i= n_0\circ \iota$, and an isomorphism
$\iota^\ast n_0^\ast \, \omega_{C_h/k} \cong \omega_{C_h^+/k}(Q_1+\ldots+Q_m)$, where $Q_1,\ldots,Q_m$ denote the marked points on $C_h^+$. Thus there is a natural isomorphism of sheaves
\[  i^\ast \omega_{{\cal C}_0/\hilb{g^{+},n,m}} \:\: \cong \:\: \omega_{{\cal C}_{{g^{+},n,m}}/\hilb{g^{+},n,m}}({\cal S}_1+\ldots+{\cal S}_m) .\]
There is a natural surjection of sheaves $\omega_{{\cal C}_0 / \hilb{g^{+},n,m}} \rightarrow i_\ast i^\ast \omega_{{\cal C}_0/\hilb{g^{+},n,m}}$. After applying $n$-th tensor powers, and push-forward, we obtain a surjection 
\[ (u_0)_\ast(\omega_{{\cal C}_0/\hilb{g^{+},n,m}})^{\otimes n} \:\: \longrightarrow \:\: u_\ast^+ i^\ast (\omega_{{\cal C}_0/\hilb{g^{+},n,m}})^{\otimes n} \]
of locally free sheaves, 
using the factorization $u^+ = u_0\circ i$. Together with the above isomorphism, and applying Grothendieck's construction of the associated projective bundle, we obtain an inclusion of projective bundles
\[ \PP u^{+}_\ast(\omega_{{\cal C}_{{g^{+},n,m}}/\hilb{g^{+},n,m}}({\cal S}_1+\ldots+{\cal S}_m))^{\otimes n} \hookrightarrow  \PP(u_0)_\ast(\omega_{{\cal C}_0/\hilb{g^{+},n,m}})^{\otimes n} .\]
Recall that $f_0 : C_0 \rightarrow \Spec(k)$ is the stable curve fixed throughout this section, and $f_0^{+} : C_0^{+}  \rightarrow \Spec(k)$ is the subcurve defined in \ref{notc+}. We fix once and for all an embedding of $C_0$ into $\PP^N$, which is equivalent to distinguishing an isomorphism 
$\PP(f_0)_\ast(\omega_{C_0/k})^{\otimes n} \cong \PP^N $, 
and we fix  an embedding of $C_0^{+}$ into $\PP^{N^{+}}$, with $N^+ := (2g^+-2+m)n-g^++1$, corresponding to an isomorphism
$ \PP(f_0^{+})_\ast(\omega_{C_0^{+}/k}(Q_1+\ldots+ Q_m))^{\otimes n} \cong \PP^{N^{+}} $, 
in such a way that the diagram
\[\begin{array}{cccccc}
& C_0^+ &\hookrightarrow & \PP^{{N^+}} & \cong & \PP(f_0^{+})_\ast(\omega_{C_0^{+}/k}(Q_1+\ldots+ Q_m))^{\otimes n}\\
(\diamond)\qquad &\downarrow &&\downarrow&&\downarrow\\[3pt]
& C_0 &\hookrightarrow & \PP^{{N}} & \cong & \PP(f_0)_\ast(\omega_{C_0/k})^{\otimes n}
\end{array}\] 
commutes. 
The universal curve $u^{+}: 
{\cal C}_{{g^{+},n,m}} \rightarrow
\hilb{g^{+},n,m}$ determines a trivialization 
\[ \begin{array}{c}
\PP u^{+}_\ast (\omega_{{\cal C}_{{g^{+},n,m}}/\hilb{g^{+},n,m}}({\cal S}_1+\ldots+{\cal S}_m))^{\otimes n} \cong \PP^{N^{+}} \times \hilb{g^{+},n,m} \hfill\\[4mm]
\hfill \cong \PP (f_0^{+})_\ast (\omega_{{C_0^{+}/k}}(Q_1+\ldots+Q_m))^{\otimes n}  \times \hilb{g^{+},n,m}.
\end{array}
\]
The tautological section of this projective bundle shall be denoted by
\[ \tau^+ \: : \quad \hilb{g^{+},n,m} \: \longrightarrow \: \PP u^{+}_\ast (\omega_{{\cal C}_{{g^{+},n,m}}/\hilb{g^{+},n,m}}({\cal S}_1+\ldots+{\cal S}_m))^{\otimes n}. \] 
If $[C^+] \in \hilb{g^{+},n,m}$ is a point representing an embedded  curve $C^+$, then by definition of $\tau^+$, the point $\tau^+([C^+]) $ is fixed under the action of the group
$ \Aut(C^+) \cong \Stab_{\PGL(N^++1)}([C^+]) $
on the fibre of the projective bundle over the point $[C^+]$. 
\end{abschnitt}

\begin{rk}\em\label{R311}\label{pgl_plus}
Via diagram $(\diamond)$ of the above construction, we may consider $\PGL(N^{+}+1)$ as a subgroup of $\PGL(N+1)$.  The group of automorphisms of $C_0$ can be identified with a subgroup of $\PGL(N+1)$ via
$ \Aut(C_0) \cong \Stab_{\PGL(N+1)}([C_0])$, and the diagram
\[\begin{array}{ccccc}
 \Aut(C_0^+) &\cong&\Stab_{\PGL(N^+ +1)}([C_0^+]) & \subset& \PGL(N^+ +1) \\[2pt]
\downarrow&&\downarrow&&\downarrow\\
\Aut(C_0) & \cong &\Stab_{\PGL(N+1)}([C_0])&\subset & \PGL(N +1)
\end{array}\]
commutes. 
 Furthermore, for any embedded curve $C^+$, which is represented by a point $[C^+]\in \hilb{g^{+},n,m}$, there is a natural  embedding of $\Aut(C^+)$ first into the group $\PGL(N^+ +1)$, and from there into $\PGL(N+1)$. 
\end{rk}

\begin{abschnitt}\em\label{const_part_2}
{\em Construction step 2.} 
As in \ref{const_part_1}, there is also a natural surjection  of the sheaves
$ j^\ast \omega_{{\cal C}_0/\hilb{g^+,n,m}} \rightarrow  \omega_{C_0^-\times\hilb{g^+,n,m}/\hilb{g^+,n,m}}({\cal S}_1+\ldots+{\cal S}_m) $ , 
where $j : C_0^-\times\hilb{g^+,n,m} \rightarrow {\cal C}_0$ denotes the embedding as a subcurve over $\hilb{g^+,n,m}$ as  constructed above. Clearly there is a trivialization
\[ \PP^{N^-}\times \hilb{g^+,n,m} \cong \PP ( (\pr_2)_\ast (\omega_{C_0^-\times\hilb{g^+,n,m}/\hilb{g^+,n,m}}({\cal S}_1+\ldots+{\cal S}_m))^{\otimes n}) \]
as a projective bundle over $\hilb{g^+,n,m}$, where $N^- = (2g^- -2+m)n-g^- +1$, provided $n$ has been chosen sufficiently large. We may assume that the corresponding  embedding of $C_0^-$ into $\PP^{N^-}$ is generic in the sense that the $m$ marked points on $C_0^-$ are contained in no subspace of dimension $m-2$. The tautological section shall be denoted by $\tau^-$.  

The above projective bundle is again in a natural way a subbundle of the  projective bundle
$\PP(u_0)_\ast(\omega_{{\cal C}_0/\hilb{g^+,n,m}})^{\otimes n} $. 
\end{abschnitt}

\begin{abschnitt}\em \label{const_part_3}
{\em Construction step 3.} 
Combining the constructions of \ref{const_part_1} and \ref{const_part_2}, we obtain two subbundles of $\PP(u_0)_\ast(\omega_{{\cal C}_0/\hilb{g^+,n,m}})^{\otimes n} $, together with  a surjective morphism $\pi$:
 \[ 
\begin{array}{ccc}
\PP u^{+}_\ast (\omega_{{\cal C}_{{g^{+},n,m}}/\hilb{g^{+},n,m}}(\sum_{i=1}^m{\cal S}_i))^{\otimes n} \\[3mm]
\oplus &\stackrel{\pi}{\longrightarrow} & \PP(u_0)_\ast( \omega_{{\cal C}_0/\hilb{g^+,n,m}})^{\otimes n}  \\
\PP (\pr_2)_\ast (\omega_{C_0^-\times\hilb{g^+,n,m}/\hilb{g^+,n,m}}(\sum_{i=1}^m{\cal S}_i))^{\otimes n}
\end{array}  \]
of projective bundles over $\hilb{g^+,n,m}$. Note that the intersection of the two subbundles can be described as the projectivization  of the cokernel of an injective  homomorphism 
\[ 
\begin{array}{ccc}
&& u^{+}_\ast (\omega_{{\cal C}_{{g^{+},n,m}}/\hilb{g^{+},n,m}}({\cal S}_1+\ldots+{\cal S}_m))^{\otimes n} \\[3mm]
 (u_0)_\ast( \omega_{{\cal C}_0/\hilb{g^+,n,m}})^{\otimes n}&\rightarrow & \oplus   \\
&& (\pr_2)_\ast (\omega_{C_0^-\times\hilb{g^+,n,m}/\hilb{g^+,n,m}}({\cal S}_1+\ldots+{\cal S}_m))^{\otimes n}
\end{array}  \]
of locally free sheaves over the base $\hilb{g^+,n,m}$. 

Thus, for a fibre of $u_0: {\cal C}_0\rightarrow  \hilb{g^{+},n,m}$, the subcurve $C_0^-$ is embedded into a subspace $\PP^{N^-} \subset \PP^N$, while the closure of its complement is embedded into a subspace $\PP^{N^+} \subset \PP^N$, and the intersection of $\PP^{N^-}$ and $\PP^{N^+}$ is a subspace of dimension $m-1$. This is just the subspace spanned by the $m$ marked points of $C_0^-$ embedded into $\PP^N$. 

By construction, for the two sections of the subbundles holds
$ \pi \circ \tau^+ \: = \: \pi\circ \tau^- $, considered 
as a  section from $\hilb{g^+,n,m}$ into the bundle $\PP(u_0)_\ast(\omega_{{\cal C}_0/\hilb{g^+,n,m}})^{\otimes n} $. In \ref{const_part_4} we will use this section to construct a trivialization of the ambient $\PGL(N+1)$-bundle, extending the trivializations of the subbundles. 
\end{abschnitt}

\begin{rk}\em \label{R311b}
In exactly the same way as in remark \ref{pgl_plus}, the construction of paragraph \ref{const_part_2} induces an embedding of $\PGL(N^- +1)$ into $\PGL(N+1)$. In particular, the group $ \Aut(C_0^-;\CQ)$ can be considered as a subgroup of $\PGL(N+1)$, too. 
\end{rk}

\begin{notation}\em\label{def_split_h}
We define the nonempty reduced subscheme
\[ \hx \: \subset \: \hilb{g^{+},n,m}  \]
of the Hilbert scheme as  the locus of points $[C^+]\in\hgprime$ parametrizing curves  $f^+:C^+\rightarrow \Spec(k)$ with the following property: if $f: C\rightarrow \Spec(k)$ denotes the stable curve of genus $g$, which is obtained by glueing $C^+$ with $C_0^-$ in the $m$ marked points, then its group of automorphisms splits naturally as
$ \Aut(C) = \Aut_\Gamma(C^+;\calq)\times \Aut(C_0^-) $. It follows from  remark \ref{newsplitrk} that $\hx$ is a dense subscheme of $\hgprime$. 

Let $\mx$ denote the GIT-quotient of $\hx$ with respect to the natural $\PGL(N^+ +1)$-action. For $\Gamma = \Gamma(C_0^-;\calq)$ we define $\hxgamma := \hx/\Gamma$ and $\mxgamma := \mx/\Gamma$. 
\end{notation} 

\begin{abschnitt}\em\label{const_part_4}
{\em Construction step 4.}
Let $\pi$ denote the morphism of projective bundles over $\hilb{g^{+},n,m} $ as in \ref{const_part_3}. The restriction of $\pi$ to $\hx$, composed with the tautological sections $\tau^+$ and $\tau^-$ of the subbundles, respectively,  defines in both cases the same  section 
$ \tau  :   \hx  \: \longrightarrow \: \PP(u_0)_\ast( \omega_{{\cal C}_0/\hilb{g^+,n,m}})^{\otimes n}|\hx$, 
compare \ref{const_part_3}. 
By \ref{const_part_1} and \ref{const_part_2}, there exist trivializations
\[ \theta^+: \quad \PP u^{+}_\ast (\omega_{{\cal C}_{{g^{+},n,m}}/\hilb{g^{+},n,m}}(\sum_{i=1}^m{\cal S}_i))^{\otimes n}|\hx \: \rightarrow \: \PP^{N^+} \times  \hx \]
and 
\[ \theta^-:  \PP (\pr_2)_\ast (\omega_{C_0^-\times\hilb{g^+,n,m}/\hilb{g^+,n,m}}(\sum_{i=1}^m{\cal S}_i))^{\otimes n}|\hx  \rightarrow \PP^{N^+} \times  \hx \]
as $\PGL(N^+ + 1)$-bundles and $\PGL(N^- + 1)$-bundles, respectively. We want to extend both of these to a common trivialization of the ambient $\PGL(N+1)$-bundle $\PP(u_0)_\ast( \omega_{{\cal C}_0/\hilb{g^+,n,m}})^{\otimes n}|\hx$. If $e$ is a point in the fibre of $\PP(u_0)_\ast( \omega_{{\cal C}_0/\hilb{g^+,n,m}})^{\otimes n}|\hx$ over a point $[C^+]\in \hx$, representing an embedded $m$-pointed stable curve $C^+$, then $e$ can be written as
$ e = \gamma \cdot \tau([C^+]) $
for some $\gamma \in \PGL(N+1)$. We now define 
\[ \theta : \quad  \PP(u_0)_\ast( \omega_{{\cal C}_0/\hilb{g^+,n,m}})^{\otimes n}|\hx \: \longrightarrow \: \PP^{N} \times  \hx\]
by
$ \theta(e) := \gamma\cdot \theta^+(\tau^+([C^+])) $, 
using the fact that the point $\tau([C^+]) = \pi\circ \tau^+([C^+])$ is contained in the subbundle $\PP u^{+}_\ast (\omega_{{\cal C}_{{g^{+},n,m}}/\hilb{g^{+},n,m}}(\sum_{i=1}^m{\cal S}_i))^{\otimes n}$, and the fixed embedding of the subspace  $\PP^{N^+} \subset \PP^N$ from \ref{const_part_1}. 
Note that in general $\gamma$ is not uniquely determined by the point $e$. To show that the above map is indeed well-defined, it suffices to show that all elements of $\PGL(N+1)$, which stabilize the point $\tau([C^+])$ also stabilize the point $\theta^+(\tau^+([C^+]))$. By construction, the elements stabilizing the point $\tau([C^+])$ in the fibre are exactly the elements of the group
$ \Aut(C) \cong \Stab_{\PGL(N+1)}([C])$, 
where $C$ is the embedded curve in $\PP^N$  obtained by glueing $C^+$ and $C_0^-$ in their marked points.  Recall that $\tau = \pi\circ\tau^+ = \pi\circ\tau^- $, and thus $\theta^+(\tau^+([C^+]))=\theta^-(\tau^-([C^+])) \in \PP^N$, with respect to the fixed  embeddings of \ref{const_part_1} and \ref{const_part_2}. 

By the definition of the subscheme $\hx$ we have a splitting of $ \Aut(C) $ as a product $ \Aut_\Gamma(C^+;\calq)\times \Aut(C_0^-)$. Because $\tau^+$ and $\tau^-$ are the tautological sections of the subbundles, the points $\theta^+(\tau^+([C^+]))$ and  $\theta^-( \tau^-([C^+]))$ are fixed by  $ \Aut(C^+) $ and $ \Aut(C_0^-)$, respectively. Since the process of glueing makes the ordering of the $m$ glueing points contained in $\calq$ disappear, the point $\theta^+ (\tau^+([C^+]))$ is also fixed under the action of $ \Aut_\Gamma(C^+;\calq) $.
Therefore, the point $\theta^+(\tau^+([C^+])) = \theta^-( 
\tau^-([C^+]))$ is fixed unter the action of  $  \Aut(C) $, and we are done. 

We have thus constructed a trivialization of the restricted projective bundle
\[ \begin{array}{rcl}
\PP(u_0)_\ast(\omega_{{\cal C}_0/\hilb{g^{+},n,m}})^{\otimes n} \: | \: \hx
  & \cong &\PP(f_0)_\ast(\omega_{{C}_0/k})^{\otimes n}
\times 
\hx\\[3mm]
&\cong&
\PP^N \times 
\hx,
\end{array} \]
which is fibrewise compatible with the natural inclusions of both projective bundles 
$ \PP (f_0^{+})_\ast (\omega_{{C_0^{+}/k}}(Q_1+\ldots+Q_m))^{\otimes n} $ and
$ \PP (f_0^{-})_\ast (\omega_{{C_0^{-}/k}}(Q_1+\ldots+Q_m))^{\otimes n}$ into   $ \PP(f_0)_\ast(\omega_{{C}_0/k})^{\otimes n}$, as constructed 
above. Globally, there is a commutative diagram of projective bundles over $\hx$: 
\[  \diagram 
\PP^{{N^{+}}} \times \hx \dto|<\hole|<<\ahook  \quad \cong & 
\PP u^{+}_\ast (\omega_{{\cal C}_{{g^{+},n,m}}/\hilb{g^{+},n,m}}(\sum_{i=1}^m{\cal S}_i))^{\otimes n} |\hx\dto|<\hole|<<\ahook\\
\PP^{{N}} \times \hx \quad \: \: \cong  & 
\PP(u_0)_\ast(\omega_{{\cal C}_0/\hilb{g^{+},n,m}})^{\otimes n}|\hx\\
\PP^{{N^{-}}} \times \hx \uto|<\hole|<<\ahook \quad\cong &
 \PP ( (\pr_2)_\ast (\omega_{C_0^-\times\hx/\hx}(\sum_{i=1}^m{\cal S}_i))^{\otimes n}).\uto|<\hole|<<\ahook 
\enddiagram \]
By the universal property of the Hilbert scheme $ \hilb{g^{+},n,m}$, the trivialization of the projective bundle $\PP(f_0)_\ast(\omega_{{C}_0/k})^{\otimes n}$ restricted to $ \hx$ 
defines a morphism
$ \thx :  \hx \rightarrow \hilb{g,n,0} $, 
which in general is not an embedding. By construction, this morphism is equivariant with respect to the action of the group $\PGL(N^+ +1)$, considered as a subgroup of $\PGL(N+1)$ as in remark \ref{pgl_plus}.

We summarize our construction in the following proposition. Using the universal property of Hilbert schemes, it is formulated geometrically in terms of embedded curves rather than those of projective bundles, while the content is of course equivalent to our construction.
\end{abschnitt}

\begin{prop}\label{const_prop}
Let ${\cal C}^\times_{{g^{+},n,m}} \rightarrow
\hx$ be the  the restriction of the universal embedded $m$-pointed stable curve  on the  Hilbert scheme $\hilb{g^{+},n,m}$ to $\hx$, and let $u_0^\times : {\cal
C}_0^\times \rightarrow  \hx$ denote the stable curve of genus $g$ obtained by glueing it along the $m$ given sections to the
trivial $m$-pointed prestable curve $C_0^{-}\times \hx
\rightarrow \hx$. 

Then there exists a global  embedding of the curve $u_0^\times : {\cal
C}_0^\times \rightarrow  \hx$ into $\pr_2: \PP^N\times \hx \rightarrow \hx$, which makes the following diagram commutative:
\[ \diagram 
{\cal C}_{g^{+},n,m}^\times\rto|<\hole|<<\ahook \dto|<\hole|<<\ahook 
& {\cal C}_0 \dto|<\hole|<<\ahook
&C_0^-\times  \hx\lto|<\hole|<<\ahook \dto|<\hole|<<\ahook 
\\
\PP^{N^{+}} \times \hx \drto \rto|<\hole|<<\ahook 
&\PP^N \times \hx \dto
&\PP^{N^{-}} \times \hx \dlto \lto|<\hole|<<\ahook \\
& \hx. 
\enddiagram \]
In particular, there exists a $\PGL(N^+ +1)$-equivariant morphism
\[ \thx : \quad \hx \rightarrow \hilb{g,n,0} . \]
\end{prop}

\proof The construction was given in four steps in  the  paragraphs \ref{const_part_1}, \ref{const_part_2}, \ref{const_part_3},  and \ref{const_part_4}.
\ebew

\begin{prop}\label{P318}
The  morphism $\thx: \hx\rightarrow \hilb{g,n,0}$ from proposition \ref{const_prop} is invariant under the action of $\, \Gamma(C_0^-;\CQ)$, 
\end{prop}

\proof
Recall that $\Gamma(C_0^-;\CQ)$ acts freely on $\hilb{g^+,n,m}$. The projective bundle $\PP u^{+}_\ast (\omega_{{\cal C}_{{g^{+},n,m}}/\hilb{g^{+},n,m}}({\cal S}_1+\ldots+{\cal S}_m))^{\otimes n}$ is invariant under the action of $\Sigma_m$, and hence in particular under the action of $\Gamma(C_0^-;\CQ)$. This is just saying that the embedding of a fibre of $u^{+} : {\cal C}_{{g^{+},n,m}} \rightarrow \hilb{g^{+},n,m}$  into $\PP^{N^{+}}$ does not depend on the order of the marked points. By construction \ref{const_part_1}, the embedding of this projective bundle into $\PP (u_0)_\ast (\omega_{{\cal C}_{0}/\hilb{g^{+},n,m}})^{\otimes n}$ is also independent of the ordering of the marked points. Thus for any point $[C']\in \hx$ representing an embedded curve $C'\subset \PP^{N^+}$, the embedding of $C'$ into $\PP^N$ is the same as the embedding corresponding to the point $\gamma([C'])$, for any  $\gamma\in \Gamma(C_0^-;\CQ)$. So it remains to consider the complement of $C'$ in the curve $C$ representing the point $\thx([C'])$. The closure of the complement of $C'$ in $C$ is just the distinguished embedding of the curve $C_0^-$, by construction of the morphism $\thx$. By definition of $\Gamma(C_0^-;\CQ)$, and by the choice of the embedding of $\Aut(C_0^-;\CQ)$ into $\PGL(N+1)$ this complement is also invariant under the action of $\Gamma(C_0^-;\CQ)$. Therefore the whole curve $C\subset \PP^N$, and hence the point $\thx([C']) \in \hilb{g,n,0}$,  is fixed under the action of $\Gamma(C_0^-;\CQ)$. 
\ebew

\begin{cor}\label{C318}
The  morphism $\thx: \hx\rightarrow \hilb{g,n,0}$ is equivariant with respect to the action of the group $\PGL(N^+ +1)\times \Gamma(C_0^-;\CQ)$.
\end{cor}

\proof The equivariance with respect to the individual groups follows from proposition \ref{const_prop} and proposition \ref{P318}. Note that the actions of the groups $\PGL(N^+ +1)$ and $\Gamma(C_0^-;\CQ)$ on $\hx$ commute, since the embedding of an $m$-pointed stable curve into $\PP^{N^+}$ is independent of the ordering of its marked points. 
\ebew 

\begin{rk}\em\label{threefoldproduct}
We can view the group $\Aut(C_0^{-})\cong \Stab_{\PGL(N^-+1)}([C_0^-])$ as a subgroup of $\PGL(N^- +1)$, and thus as a subgroup of $\PGL(N+1)$. The inclusion is determined by the fixed embedding of $\PP^{N^-}$ into $\PP^N$. By definition, elements of $\Aut(C_0^-)$ fix the $m$ marked points of $C_0^-$, which span the intersection $\PP^{N^-} \cap \PP^{N^+}$ in $\PP^N$ by remark \ref{const_part_3}. Thus such an element acts as the identity on the intersection, and therefore elements of $\Aut(C_0^-)$ and elements of $\PGL(N^+ +1)$ commute inside $\PGL(N+1)$. In particular, we can view $\PGL(N^{+}+1) \times \Aut(C_0^{-}) $ as a subgroup of $\PGL(N+1)$. 

Furthermore, the action of $\Aut(C_0^-)$ on $ \hilb{g,n,0}$ commutes with the action of $\gammaco$, since the embedding of $C_0^-$ into $\PP^N$ is independent of the labels of the $m$ marked points. We can thus strenghten corollary \ref{C318} as follows.
\end{rk}

\begin{lemma}\label{L42}
The morphism $\thx: \hx \rightarrow \hilb{g,n,0}$ is
equivariant with respect to the action of $\, \PGL(N^{+}+1) \times \gammaco \times \Aut(C_0^{-})$, where $\Aut(C_0^{-})$ acts trivially on $\hx$.
\end{lemma}

\proof
It only remains to show that the image of $\hx$ under $\thx$ is pointwise fixed under the action  of $\Aut(C_0^{-})$. Recall that  $\thx$ is constructed  via  compatible trivializations as in the following diagram:
\[ \diagram
\PP u^{+}_\ast(\omega_{{{\cal
C}}_{g^{+},n,m}/\hilb{g^{+},n,m}}(\sum_{i=1}^m{\cal S}_i ))^{\otimes n}  \rto|<\hole|<<\ahook \dto^{\cong}& \PP(u_0)_\ast(\omega_{{\cal
C}_0/\hilb{g^{+},n,m}})^{\otimes n} \dto_{\cong}\\
\PP(f_0^{+})_\ast (\omega_{C_0^{+}/k}(\sum_{i=1}^m Q_i))^{\otimes
n}  \times \hilb{g^{+},n,m} \rto|<\hole|<<\ahook &
\PP(f_0)_\ast(\omega_{C_0/k})^{\otimes n} \times
\hilb{g^{+},n,m},
\enddiagram \]
restricted to the subscheme $\hx$. 
By definition of the inclusion of $\Aut(C_0)$ into $\PGL(N+1)$,  the embedding of $C_0$ into $\PP^N$ is fixed under the action of $\Aut(C_0)$, and hence in particular fixed under the action of $\Aut(C_0^{-})$. Therefore the embedding of the projective bundles is invariant under the action of $\Aut(C_0^{-})$, and thus $\thx$ is also invariant with respect    to the action of $\Aut(C_0^{-})$. 
\ebew

\begin{prop}\label{312}
Let $\, \Gamma := \gammaco$. 
The  morphism $\thx$
induces an embedding  
\[ \Theta : \quad  \hxgamma   \: \rightarrow \:
\hilb{g,n,0}  \]
and an
isomorphism of schemes
$ M_{g^{+},\mgamma}^\times \cong  \dpsplit $. 
\end{prop}

 \proof
Proposition \ref{P318} implies that the morphism $\thx : \hx \rightarrow \hilb{g,n,0}$ factors through $\hxgamma = \hx/\Gamma$. 

Let $C_1$ and $C_2$ denote two embedded $m$-pointed stable curves over $\Spec(k)$, represented by two points $[C_1], [C_2]\in \hx$. Suppose that $\thx([C_1]) = \thx([C_2]) \in \hilb{g,n,0}$, and let $C$ denote the embedded curve representing this point. By the assumption on $[C_1]$ and $[C_2]$, the curve $C$ has $3g-4$ nodes, so there is a unique irreducible component on $C$, which has not the maximal number of nodes when considered as a pointed stable curve itself. This component is of course equal to the embedded curves $C_1$ and $C_2$ in $\PP^N$, considered as curves with $m$ distinguished points, but without ordered labels on these points. The closure of its complement in $C$ is isomorphic to $C_0^-$ by definition of the morphism $\thx$. Since $\thx([C_1])$ and $\thx([C_2])$ represent the same embedded curve, the order of the marked points of $C_1$ and $C_2$ can differ at most by a reordering induced by an automorphism of $C_0^-$. Hence by definition of $\Gamma(C_0^-;\CQ)$, the points $[C_1]$ and $[C_2]$ are the same modulo the action of $\Gamma(C_0^-;\CQ)$. 
This proves the first part of the proposition.

For the second part note that by corollary \ref{C318} the morphism $\thx$ induces a morphism
$   \mxgamma \rightarrow \msbar_g $.  
By the definition of $M_{g^+,m}^\times$, its image is $\dpsplit$. The morphism is just the restriction of the isomorphism from lemma \ref{decompose}. 
\ebew

\section{Decomposition of stable curves}\label{sectdecompose}

\begin{notation}\em
Let $\pi : \hilb{g,n,0} \rightarrow \msbar_g$ be the canonical
quotient morphism. We denote by 
\[ {K}(C_0;P_1) \: \subset \: \hilb{g,n,0}\]
the reduction of the preimage of $\dpo$. Furthermore, we denote  by $\kx$ the reduction of the preimage of the open subscheme $\dx$, which represents curves with splitting automorphism groups.
\end{notation}

\begin{lemma}\label{314}
The subscheme $\kx$ is equal to the reduction of the $\, \PGL(N+1)$-orbit 
$ \kx = \thx(\hx) \cdot \PGL(N+1) $
of the image of  the morphism $\thx : \hx \rightarrow \hilb{g,n,0}$ from proposition \ref{const_prop}. 
\end{lemma}

\proof
This follows immediately from the definiton of $\kx$, and the surjectivity of the morphism $M_{g^+,m}' \rightarrow \dpo$. 
\ebew

\begin{rk}\em\label{Ksmooth}
The scheme $\kpo$  is a smooth and irreducible subscheme  of $\overline{H}_{g,n,0}$. This follows from the definition of $\dpo$ as a connected  component of $M_g^{(3g-4)}$, and the fact that the boundary $\hilb{g,n,0} \setminus {H}_{g,n,0}$ is a normal crossing divisor, compare \cite{DM}. 
\end{rk}

\begin{lemma}\label{315}
Let $f : C \rightarrow S$ be a stable curve of genus $g$, with reduced and irreducible  base $S$, such that for all closed points $s\in S$ the fibre $C_s$ has exactly $3g-4$ nodes. Then the induced morphism $\theta_f : S \rightarrow \msbar_g$ factors through $\dpo$ if and only if 
there exists an \'etale  covering $u: S'\rightarrow S$, such that the pullback stable curve $f': C'\rightarrow S'$ decomposes as follows. There is an $m$-pointed prestable subcurve, which is unique up to permutations of the labels of its marked points, and isomorphic to the trivial curve $C_0^-\times S' \rightarrow S'$, together with an $m$-pointed stable curve $f^\# : C^\# \rightarrow S'$ of genus $g^+$, such that the curve $f': C'\rightarrow S'$ is obtained by glueing the two subcurves along the corresponding sections of marked points.
\end{lemma} 

\proof
Suppose that such an \'etale covering exists. Since $f':C'\rightarrow S'$ is the pullback of $f:C\rightarrow S$, it too is a stable curve of genus $g$ with exactly $3g-4$ nodes in each fibre, and the induced morphism $\theta_{f'} : S' \rightarrow \msbar_g$ factors as $\theta_{f'}= \theta_{f} \circ u$. Thus $\theta_{f}$ and $\theta_{f'}$ map to the same irreducible component $D$ of $M_g^{(3g-4)}$, since $S$ is irreducible. For a general point $s\in S'$, the closure of the complement of $C_0^-$ in  the fibre $C_s'$ is an $m$-pointed stable curve of genus $g^+$. Thus by definition $\theta_{f'}(s) \in \dpo$, and hence $D=\dpo$. 

Conversely, suppose that $\theta_f : S \rightarrow \msbar_g$ factors through $\dpo$. By lemma \ref{sectnod} there exists an \'etale covering $u:S'\rightarrow S$, such that $S'$ is connected, and the pullback $f':C'\rightarrow S'$ admits $3g-4$ sections of nodes. For a closed point $s\in S'$ there are two types of nodes on the fibre $C_s'$: those nodes, which are intersection points of different irreducible components of $C_s'$, and those which are singularities of one irreducible component. Since the curve $f':C'\rightarrow S'$ is locally a deformation of its fibre $C_s'$, and since the number of nodes is constant for all fibres, all nodes of one section must be of the same type. 

Suppose that $\sigma_1,\ldots,\sigma_k: S'\rightarrow C'$ are the sections of nodes, which are intersection points. Let $\tilde{f}: \tilde{C}\rightarrow S'$ denote the flat family, which is obtained by ``cutting'' the curve $C'$ along this $k$ sections. Formally, we may construct $\tilde{C}$ by choosing an embedding of $f': C'\rightarrow S'$ into $\PP^\ell \times S' \rightarrow S'$, for some sufficiently large $\ell$, and blowing up along the subschemes $\sigma_1(S'),\ldots,\sigma_k(S')$. The family $\tilde{f}: \tilde{C}\rightarrow S'$ is then given as the strict transform of $f': C'\rightarrow S'$, where the points of $C'$ contained in the $k$ sections are separated into pairs of points of  $\tilde{C}$ above. 

Note that $\tilde{f}: \tilde{C}\rightarrow S'$ is a pointed prestable curve over $S'$. In particular, since $S'$ is connected, for any connected component of $Z \subset \tilde{C}$ the restriction $\tilde{f}|Z : Z \rightarrow S'$  is still a pointed stable curve. Thus the genus and the number of nodes remains constant for all fibres of $\tilde{f}|Z$. For a general closed point $s \in S'$, the fibre $C_s'$ has exactly $3g-4$ nodes. Hence there exists a unique irreducible component of $C_s'$, which has not the maximal number of nodes when considered as a pointed stable curve itself. Therefore there exists a unique connected component of the fibre $\tilde{C}_s$, which has not the maximal number of nodes when considered as a pointed stable curve. This component is necessarily an $m$-pointed stable curve of genus $g^+$. This implies that there exists a unique connected component $C^\#$ of $\tilde{C}$, which is an $m$-pointed stable curve of genus $g^+$ over $S'$. By construction, this curve $f^\#: C^\#\rightarrow S'$ can be considered as a subcurve embedded in  $f': C'\rightarrow S'$. 

Let $f^-: C^-\rightarrow S'$ denote the family obtained as the closure of the complement of $f^\#: C^\#\rightarrow S'$ in  $f': C'\rightarrow S'$. Let $s\in S'$ be a closed point, such that $\theta_{f'}(s) \in \dpo$. Then the fibre $C_s'$ must contain a subcurve isomorphic to $C_0^-$ by definition. All irreducible components of $C_0^-$ have the maximal number of nodes when considered as pointed stable curves, so $C_0^-$ must be equal to the closure of the complement of $C_s^\#$ in $C_s'$. Thus, all fibres of the $m$-pointed prestable curve $f^-: C^-\rightarrow S'$ are isomorphic to $C_0^-$. In particular, after applying another \'etale covering of $S'$ if necessary, the curve $f^-: C^-\rightarrow S'$ is isomorphic to the trivial family $C_0^-\times S' \rightarrow S'$. 
\ebew

Note that although  $f^\# : C^\# \rightarrow S'$ is uniquely determined as a subscheme, being the complement of $C_0^-\times S' \rightarrow S'$, the labels of its $m$ sections of marked points depend on choices made during the process of glueing. 

\begin{cor}\label{K324}
Let $f: C \rightarrow S$ be a stable curve of genus $g$, with reduced base $S$. Suppose that
the induced morphism $\theta_f : S \rightarrow \msbar_g$ factors through
$\dpo$. Then there exists a uniquely determined decomposition of $f: C\rightarrow S$ into an $m/\gammaco$-pointed stable curve $f^+:C^+\rightarrow S$ of genus $g^+$, and an $m/\gammaco$-pointed prestable curve $f^-:C^-\rightarrow S$ of genus $g^-$, where the second curve is locally isomorphic to the trivial curve $C_0^-\times S\rightarrow S$ with respect to the \'etale topology. 
\end{cor}

\proof
Without loss of generality we may assume that $S$ is irreducible. Indeed, if decompositions exist on the irreducible components of $S$, then the uniqueness property implies that they fit together to a global decomposition. 

By lemma \ref{315} there exists an \'etale covering $u:S'\rightarrow S$, such that the pullback curve $f':C'\rightarrow S'$ decomposes into an $m$-pointed stable curve $f^\#: C^\#\rightarrow S'$, and a trivial $m$-pointed prestable curve $C_0^-\times S' \rightarrow S'$. By assumption, each fibre of $f:C\rightarrow S$ has exactly $3g-4$ nodes. Consider the pullback morphism $\overline{u}:C'\rightarrow C$. For a closed point $s\in S'$, it maps the fibre $C^\#_s$ of the subcurve $f^\#:C^\#\rightarrow S'$ isomorphically onto the unique irreducible component of the fibre $C_{u(s)}$, which has not the maximal number of nodes when considered as a pointed stable curve. In particular, the image of $C^\#$ in $C$ defines uniquely a subfamily $f^+: C^+\rightarrow S$ in $f:C\rightarrow S$. 

Analogously, the image of the trivial subcurve $C_0^-\times S'\rightarrow S'$ defines a subfamily $f^-: C^-\rightarrow S$. Since the latter is the \'etale image of an $m$-pointed prestable curve, it is by definition an $m/\gammaco$-pointed curve.

Let $s,s'\in S'$ be two closed points with $u(s)=u(s')$. If $C_{s}^\#$ and $C_{s'}^\#$ denote the two corresponding fibres of $f^\#:C^\#\rightarrow S'$, then glueing in the $m$ marked points to $C_0^-$ produces in both cases the same fibre $C_{u(s)} = C_{u(s')}$ of $f: C\rightarrow S$. Thus the labels of the respective $m$ marked points may differ only by a permutation induced by an automorphisms of $C_0^-$, i.e. by a permutaion contained  in $\gammaco$. Thus the induced curve $f^+: C^+\rightarrow S$ is an $m/\gammaco$-pointed stable curve of genus $g^+$, as defined in  \cite{Z2}. 
\ebew

The above corollary \ref{K324} deals only with stable curves $f:C\rightarrow S$, where the base $S$ is reduced. In our applications, we will need an analogous statement for arbitrary schemes $S$, and for this we need one extra assumption.

\begin{rk}\em\label{EPHI}
There is a functorial way to assign to each pointed stable curve $f: C\rightarrow S$ a pair $(\ecs,\pcs)$, consisting of a principal bundle $\prcs:\ecs\rightarrow S$, and an equivariant morphism $\pcs:\ecs \rightarrow \overline{H}_{g,n,m}$. We recall briefly the description of this construction from \cite{Z2}.  

Let $f:C\rightarrow S$ be an $m$-pointed stable curve of genus $g$, with sections of marked points $\sigma_1,\ldots,\sigma_m: S\rightarrow C$. The sections define divisors $S_1,\ldots,S_m$ on $C$. We denote by $\prcs:\ecs\rightarrow S$ the $\PGL(N+1)$-principal bundle associated to the projective bundle $P:= \PP f_\ast(\omega_{C/S}(S_1+\ldots+S_{m}))^{\otimes n}$ over $S$. The pullback of $P$ to $\ecs$ has a natural trivialization as a projective bundle. This implies that the pullback stable curve $\overline{f}:\overline{C}\rightarrow \ecs$ of the stable curve $f:C\rightarrow S$ has a natural embedding into $\PP^N\times \ecs$. By the universal property of the Hilbert scheme this is equivalent to specifying a $\PGL(N+1)$-equivariant morphism $\pcs: \ecs \rightarrow \overline{H}_{g,n,m}$. 

It is straightforward to see that this construction is functorial. The assignment is even reversible: let $u: {\cal C}_{g,n,m} \rightarrow \overline{H}_{g,n,m}$ denote the universal curve. Then the curve $f:C\rightarrow S$ can be recovered from the pair $(\ecs,\pcs)$ by taking the quotient $\phi^\ast_{C/S} ({\cal C}_{g,n,m}) /\PGL(N+1) \rightarrow \ecs/\PGL(N+1)$. 
\end{rk}

\begin{prop}\label{K324a}
Let $f: C \rightarrow S$ be a stable curve of genus $g$. Suppose that
the induced morphism $\theta_f : S \rightarrow \msbar_g$ factors through
$\dpo$. Let $(\ecs,\pcs)$ be the pair associated to $f: C \rightarrow S$ as in remark \ref{EPHI}. Suppose that $\pcs$ factors through $\kpo$. 
Then there exists a uniquely determined decomposition of $f: C\rightarrow S$ into an $m/\gammaco$-pointed stable curve $f^+:C^+\rightarrow S$ of genus $g^+$, and an $m/\gammaco$-pointed prestable curve $f^-:C^-\rightarrow S$ of genus $g^-$, where the second curve is locally isomorphic to the trivial curve $C_0^-\times S\rightarrow S$ with respect to the \'etale topology. 
\end{prop}

\proof
Let ${\cal C}_{g,n,0}\rightarrow \overline{H}_{g,n,0}$ denote the universal curve, and ${\cal C}_{K}\rightarrow \kpo$ its restricition to $\kpo$. Since $\kpo$ is reduced, corollary \ref{K324} 
shows the existence of a unique subcurve ${\cal C}_K^-$, which is locally isomorphic to the trivial curve $C_0^-\times \kpo$ as an $m/\Gamma(C_0^-;{\cal Q})$-pointed curve. Let ${\cal C}_E$ and ${\cal C}_E^-$ denote the  respective pullbacks of ${\cal C}_K$ and ${\cal C}_K^-$ to $\ecs$. By the definition of $\ecs$, one has $f: C={\cal C}_E/\PGL(N+1) \rightarrow \ecs/\PGL(N+1)=S$, and $C^- := {\cal C}_E^-/\PGL(N+1) \rightarrow S$  is the unique subcurve, which is locally isomorphic to $C_0^-\times S\rightarrow S$ as an $m/\Gamma(C_0^-;{\cal Q})$-pointed curve.
\ebew

\section{The moduli substack}\label{par4}\label{sectmainthm}
After our discussion on the level of moduli spaces, let us now define the corresponding substacks in the boundary of the moduli stack $\cmgbar$ of Deligne-Mumford stable curves of genus $g$. 

\begin{defi}\em
Let $\FDPtilde$ denote the preimage substack of the scheme $\dpo$ in $\mgstack{g}$ under the canonical map from $\mgstack{g}$ to $\msbar_g$. 
The stack $\FDPopen$ is defined as the reduced stack underlying  $\FDPtilde$. We define the substack $\FDPcopen$ as the reduction of the preimage of the open and dense  subscheme $\dx$.
\end{defi}

\begin{rk}\em
For any scheme $S$, the set of objects of the fibre category $\FDPtilde(S)$ is equal to the set of those stable curves $f:C\rightarrow S$ of genus $g$, such that the induced morphism $\theta_f: S\rightarrow \msbar_g$ factors through the subscheme $\dpo$. It follows from proposition \ref{K324a} that  $\FDPopen$ is a moduli stack for curves $f:C\rightarrow \Spec(k)$, which contain a subcurve isomorphic to $C_0^-$, where the marked points on $C_0^-$ correspond to nodes on $C$, and that $\dpo$ is its moduli space. 
\end{rk}

\begin{lemma}\label{140}
There are isomorphisms of stacks
\[ \FDPopen \: \cong \: \left[ \, {K}(C_0;P_1) / \PGL(N+1) \right],\]
and
\[ \FDPcopen \: \cong \: \left[ \, \kx / \PGL(N+1) \right].\]
In particular, both stacks are smooth and irreducible, and of dimension one.
\end{lemma}

\proof 
It was shown by Edidin \cite{Ed} that the moduli stack of stable curves of genus $g$ can be constructed as a quotient stack $\cmgbar \cong [\hilb{g,n,0}/\PGL(N+1)]$. Using this, the descriptions of $\FDPopen$ and $\FDPcopen$ follow immediately from the definitions by standard arguments on Cartesian products and quotient stacks. Since by remark \ref{Ksmooth} the scheme ${K}(C_0;P_1)$ is  smooth and irreducible as an atlas of $\FDPopen $, so is $\FDPopen$ itself.
\ebew

\begin{rk}\em
From lemma \ref{140} the following characterization of $\cdpo$ can be derived. For any scheme $S$, the fibre category $\cdpo(S)$ has as its objects such stable curves $f:C\rightarrow S$ of genus $g$, where  the induced morphism $\theta_f : S\rightarrow \mgbar$ factors through $\dpo$, and the morphism $\pcs: \ecs \rightarrow \hilb{g,n,0}$ from remark \ref{EPHI} factors through ${K}(C_0;P_1)$. 
\end{rk}

\begin{notation}\label{not004}\em
Let $\cmgprime$ denote the substack of $\cmsbar_{g^+,m}$ which is the reduction of the preimage of $\mgprime$. Analogously, let $\cmgprimequot$ denote the substack of $\cmgprimequot$ which is the reduction of the preimage of $\mgprimequot$. As in \ref{mprimed}, there is a natural morphism $\theta: \cmgprime \rightarrow \cmgbar$, which is given by glueing the trivial curve with fibre $C_0^-$ to any $m$-pointed stable curve. 

However, unlike the situation in the case of moduli spaces as described in lemma \ref{decompose}, there is no induced isomorphism between the moduli stacks $\cmgprimequot$ and $\cdpo$, where $\Gamma = \Gamma(C_0^-;\calq)$. It turns out that in order to relate the two stacks, we need to consider a finite quotient of  $\cmgprimequot$, see remark \ref{Rk46}. Even then, there is no global morphism from this quotient to $\cdpo$, so we must furthermore restrict ourselves to the following  open and dense substack.

Let us denote by  $\cmgammax$ the substack of $\cmgbarquot$, which is the reduction of the preimage of $\mgammax$. Analogously we define $\cmgx$ in $\cmsbar_{g^+,m}$.
\end{notation}

By lemma \ref{L42} and proposition \ref{312}, there exists a $\PGL(N^{+}+1) \times \Aut(C_0^{-})$-equivariant morphism 
$ \Theta:  \hgammax  \rightarrow \hilb{g,n,0}$, 
where $\Aut(C_0^{-})$ acts trivially on $\hgammax$. The following proposition provides a stack version of the induced morphism $\mgammax \rightarrow \dx$.

\begin{prop}\label{MainProp}
There is a morphism of Deligne-Mumford stacks
\[ \Lambda : \quad  \left[\hgammax \, / \, \PGL(N^{+}+1) \times \Aut(C_0^{-}) \right] \: \: \longrightarrow  \: \:  \FDPcopen \]
where $\Aut(C_0^{-})$ acts trivially on $\hgammax$, and $\Gamma:= \gammaco$. 
\end{prop}

As abbreviations, we will from now on use $A:= \Aut(C_0^{-})$, $\Gamma := \Gamma(C_0^-;\CQ)$  and $P := \PGL(N^{+}+1)$. We will also denote the subschemes $\kx$ and $\dx$ simply by $\kpx$ and $\dpx$, respectively.

\begin{rk}\em\label{Rk46}
In \cite{Z2}, the moduli stack of $\mgamma$-pointed stable curves has been constructed as a quotient stack $\cmgbarquot \cong [ \, \hgbarquot/P\,]$. Thus there is also an isomorphism $\cmgammax  \cong [\, \hgammax/P\,]$. By composing the morphism $\Lambda$ of proposition \ref{MainProp} with the canonical quotient morphism $\cmgammax \rightarrow \cmgammax / A = [ \, \hgammax/P \times A\,]$, we obtain a morphism of stacks 
$ \cmgammax  \rightarrow  \cmgbar $,
which factors through $\cdx$, and which induces the natural morphism $\mx \rightarrow \dx$ on the corresponding  moduli spaces. 
\end{rk}

\proofof{proposition \ref{MainProp}}
Viewing stacks as categories fibred in groupoids over the category of
 schemes, the morphism $\Lambda$  will be given as a functor respecting the fibrations.
Let a scheme $S$ be fixed. Recall that an object of 
$[ \, \hgammax  /  P\times A ](S)$ is a triple $(E',p',\phi')$, where $p':E'\rightarrow S$ is a principal $P\times A$-bundle, together with a $P\times\ A$-equivariant morphism $\phi':E' \rightarrow \hgammax$. Note that this data is equivalent to specifying a stable curve $f:C\rightarrow S$ as an object of $\cmgammax$, with $\prcs = p'$ and $\pcs = \phi'$, compare remark \ref{EPHI}. 

Viewing $P\times A $ as a subgroup of $\PGL(N+1)$ as in remark \ref{threefoldproduct}, 
we can define an extension $E$ of $E'$ by
\[ E := (E' \times \PGL(N+1) ) / (P\times A) \]
as a principal $\PGL(N+1)$-bundle $p:E\rightarrow S$ over $S$. By lemma \ref{314} we have $\kpx = \thx(\hx)\cdot\PGL(N+1) = \Theta(\hgammax)\cdot\PGL(N+1)$. Therefore we can define a morphism
\[ \begin{array}{cccc}
\phi : & E& \rightarrow & \kpx\\
& e & \mapsto& \Theta(\phi'(e'))\cdot \gamma, 
\end{array}\]
where  $e\in E$ is represented by a pair $(e',\gamma)\in E'\times \PGL(N+1)$.
This map is indeed well-defined, as $\Theta \circ\phi'$ is equivariant with respect to $P\times  A$.  For all $g \in \PGL(N+1)$ we have 
$ \phi(eg) = \phi(e) g$
for all $e\in E$, so $\phi$ is $\PGL(N+1)$-equivariant. Hence 
$ (E,p,\phi) \in \left[ \, \kpx / \PGL(N+1) \right](S) \cong \FDPcopen(S)$, 
where the last isomorphism holds by lemma \ref{140}.

It is easy to see that this construction of extending principal bundles is functorial. In this way we obtain for each  scheme $S$   a functor between  the fibre categories 
$ \Lambda_S :   [\hgammax \, / \, P\times A ](S) \rightarrow  \FDPcopen(S)$.  This in turn defines the desired  morphism of stacks
$ \Lambda : [\hgammax \, / \, P\times A ]  \rightarrow  \FDPcopen $
in the usual way. 
\ebew

We now arrive at our central result on the one-dimensional boundary strata of the moduli stack of Deligne-Mumford stable curves.

\begin{thm}\label{MainThm}
There is an isomorphism of Deligne-Mumford stacks
\[ \left[ \hgammax  \, / \, \PGL(N^{+}+1)\times \Aut(C_0^{-}) \right] \: \: \cong \: \:  \cdx \]
where $\Aut(C_0^{-})$ acts trivially on $\hgammax$, and $\, \Gamma:= \gammaco$.  
\end{thm}

\begin{rk}\em\label{rk00411}
$(i)$ By the construction given in \cite{Z2}, there is an isomorphism of stacks $\cmgammax \cong [\hgammax  \, / \, \PGL(N^{+}+1)] $. Therefore theorem \ref{MainThm} implies that generically the moduli substack $\cdpo$ is isomorphic to the stack quotient $\cmgprimequot / \Aut(C_0^-)$. Note that since the action of $\Aut(C_0^{-})$ on $\hgprimequot$ is trivial, this quotient is isomorphic to the product of the moduli stack $\cmgprimequot$ with the classifying stack $\mbox{B}\Aut(C_0^-)$ of the group $\Aut(C_0^-)$. Im particular, even though there is an isomorphism $\dx \cong \mgammax $ of the respective moduli spaces, the stacks $\cdx$ and $\cmgprimequot$ are in general not isomorphic. 

$(ii)$ Theorem \ref{MainThm} is in fact the strongest result possible. Indeed, consider  a stable curve $f: C\rightarrow \Spec(k)$, which is an object of $\cdpo(\Spec(k))$, but not of $\cdx(\Spec(k))$. An extension of the isomorphism must fail, since the splitting of the stack into a product would imply a corresponding splitting of the automorphism group $\Aut(C)$, as explained in \cite{Z1}, and thus contradict the definition of $\cdx$. 

However, theorem \ref{MainThm} can be extended to a certain degree to the closure of $\FDPopen$. Recall that $\FDPopen$ is defined by a connected component of the boundary stratum $M_g^{(3g-4)}$ parametrizing stable curves with exactly $3g-4$ nodes. So its closure contains stable curves with at least $g-4$ nodes. 
For a full discussion and the construction of this extension, see \cite{Z1}.  
\end{rk}

\proofof{theorem \ref{MainThm}}
The proof is by explicit construction of an inverse morphism $\Xi$ to the morphism $\Lambda$ from  proposition
\ref{MainProp}  in such a way, that the composition $\Lambda \circ \Xi$ is isomorphic to, and the composition $\Xi\circ \Lambda$ is equal to the respective identity functor. As before, we are using the abbreviations $A:= \Aut(C_0^{-})$, $\Gamma := \Gamma(C_0^-;\CQ)$  and $P := \PGL(N^{+}+1)$, as well as $\kpx:=\kx$ and $\dpx:= \dx$. 

We need to construct for each object of $ \cdx(S)$, for each scheme $S$, its image in $ [ \hgammax  \, / \, \PGL(N^{+}+1)\times \Aut(C_0^{-}) ](S)$ under $\Xi$. By lemma \ref{140} we may assume that such an object is given by a triple $(E,p,\phi)$, where  $p:E \rightarrow S$ is a principal $\PGL(N+1)$-bundle, together with a $\PGL(N+1)$-equivariant morphism $\phi: E \rightarrow \kpx$. In a first step we will construct from this data a principal $P\times A$-bundle $p^+: E^+\rightarrow S$. 

$(i)$ Consider the restriction $ {\cal C}_{g,n,0} \rightarrow \kpx  $ of the universal stable curve over $\hilb{g,n,0}$ to $\kpx$, 
 and denote its pullback to $E$ via $\phi$ by 
$  {\cal C}_E \rightarrow E $. 
Taking quotients by the action of $\PGL(N+1)$  induces a morphism
\[ f_S : \quad {\cal C}_S := {\cal C}_E / \PGL(N+1) \rightarrow E / \PGL(N+1) = S ,\]
which is a stable curve of genus $g$ over $S$. The induced morphism  $\theta_f : S \rightarrow \msbar_g$ factors through $\dpx$. Therefore,  by proposition \ref{K324a},  there exists a unique subcurve $f_S^-: {\cal C}_S^- \rightarrow S$ of $f_S: {\cal C}_S\rightarrow S$, and an \'etale covering $u: S'\rightarrow S$, such the pullback $f_{S'}^-:{\cal C}_{S'}^- \rightarrow S'$ of $f_S^-$  is isomorphic to the trivial $m$-pointed stable curve $C_0^-\times S' \rightarrow S'$. Let $f_{S'}:{\cal C}_{S'} \rightarrow S'$ denote the pullback of $f_S$. 

Recall that in general $C_0^{-}$ is not a stable curve, but only prestable. Suppose that $C_0^-$  decomposes into connected components
$C_1^{-},\ldots,C_r^{-}$, where $C_i^-$ is an $m_i$-pointed
stable curve of some genus $g_i^{-}$ for $i=1,\ldots,r$. This induces a decompsiton of $f_{S'}^-:{\cal C}_{S'}^- \rightarrow S'$ into subcurves $f_{S',i}^{-} : {\cal
C}_{S',i}^{-} \rightarrow S'$. Let $S_{1,i}',\ldots,S_{m_i,i}'$ denote the respective sections of the $m_i$ marked points on ${\cal C}_{S',i}^{-}$. As in construction \ref{const_part_1} there is a natural inclusion of projective bundles
$ \PP( f_{S',i}^{-})_\ast(\omega_{{\cal C}_{S',i}^{-}/S'}(S_{1,i}'+\ldots+S_{m_i,i}'))^{\otimes n} \: \hookrightarrow \:
\PP(f_{{S'}})_\ast(\omega_{{\cal C}_{S'}/S'})^{\otimes n}$, 
where we may assume that $n$ has been chosen large enough a priory. 
Put  $N_i^{-} := (2g_i^{-}-2+m_i)n-g_i^{-}+1$. There is a principal
$\PGL(N^{-}_i+1)$-bundle $p^{-}_{i} : E_{i}^{-}
\rightarrow S'$ associated to the projective
bundle $\PP( f_{S',i}^{-})_\ast(\omega_{{\cal
C}_{S',i}^{-}/S'}(S_{1,i}'+\ldots+S_{m_i,i}'))^{\otimes
n}$. 

This bundle is naturally a subbundle of the principal $\PGL(N+1)$-bundle $p_{S'}: E_{S'} \rightarrow S'$ associated to the projective bundle
$\PP(f_{{S'}})_\ast (\omega_{{\cal C}_{S'}/S'})^{\otimes n}$, which is equal to the pullback of the principal bundle $p : E\rightarrow S$ from above. Indeed, there is a natural inclusion  of $\PP( f_{S',i}^{-})_\ast(\omega_{{\cal
C}_{S',i}^{-}/S'}(S_{1,i}'+\ldots+S_{m_i,i}'))^{\otimes
n}$ into $ \PP(f_{{S'}})_\ast (\omega_{{\cal C}_{S'}/S'})^{\otimes n}$, compare \ref{const_part_1}. Viewing $A$ as a subgroup of $\PGL(N^{-}_i+1)$, this induces an embedding of the respective quotients with respect to the action of $A$, and therefore an inclusion $ E_{i}^{-}/A  \hookrightarrow E_{S'}/A$. Recall from \ref{not31} that $A$ is a normal subgroup of $\Aut(C_0^-;\calq)$, with quotient $\Aut(C_0^-;\calq)/A = \Gamma$. Since the embedding of the projective bundles is invariant with respect to reorderings of the labels of the marked points, i.e. with respect to the action of $\Gamma$, there is thus even an induced embedding
\[ E_{i}^-/\Aut(C_0;\calq) \:\: \hookrightarrow \:\: E_{S'}  / A .\] 
For $i=1,\ldots, r$ we define  $f^{-}_{E_i}: {\cal C}^{-}_{E_i} \rightarrow E_i^{-}$ as the pullback of the stable curve $f_{S',i}^{-} : {\cal
C}_{S',i}^{-} \rightarrow S'$ along $p_i^-$. Each is an $m_i$-pointed stable curve of genus $g^{-}_i$ over
$E_{i}^{-}$, with sections of marked points $\tilde{S}_{1,i}',\ldots,\tilde{S}_{m_i,i}'$. The projective bundle
$ \PP(f_{E_i}^{-})_\ast (\omega_{{\cal C}^{-}_{E_i}/{E_i}^{-}}(\tilde{S}_{1,i}'+\ldots+\tilde{S}_{m_i,i}'))^{\otimes n} $ has a natural trivialization, compare also remark \ref{EPHI}.  
Hence by the universal property of the Hilbert scheme there is an induced morphism $\theta_{i}^{-}:  E_{i}^{-} \rightarrow  \hilb{g^{-}_i,n,m_i}$. 
which is $\PGL(N^{-}_i+1)$-equivariant.

$(ii)$ Now we want to put these constructions on the connected components ${\cal C}^-_{S',i}$  together for the whole 
prestable curce  $f^-_{S'}:{\cal C}_{S'}^-\rightarrow S'$. We define $P^{-} := \PGL(N^{-}_1+1) \times \ldots \times \PGL(N^{-}_r +1)$, and $\hilb{g^{-},n,m}^{-} := \hilb{g^{-}_1,n,m_1}\times \ldots
\times  \hilb{g^{-}_r,n,m_r}$. We obtain a  principal $P^{-}$-bundle
$ p^{-}_{S'} :  E^{-}_{S'} := E_{i}^{-}\times \ldots\times E_{r}^{-} \rightarrow  S'$, together with an induced  $P^{-}$-equivariant morphism
$ \theta_{S'}^-: E_{S'}^- \rightarrow \hilb{g^{-},n,m}^{-}$. 
 
By construction, the induced morphism   $ {\cal C}^{-}_{E_1}\times \ldots\times{\cal C}^{-}_{E_r} \: / \: P^- \rightarrow E_{S'}^-/ P^-$ on the quotients is just the  $m$-pointed prestable curve $f_{S'}^-: {\cal C}_{S'}^- \rightarrow S'$ of genus $g^{-}$ over $S'$. Hence for the  induced  morphism
$ \theta':  S'  \rightarrow  \msbar_{{g^{-},m}}^{-} $ holds $\theta'\circ p_{S'}^{-} = \pi \circ \theta_{S'}^{-}$, where $\pi : \hilb{g^{-},n,m}^{-} \rightarrow \msbar_{{g^{-},m}}^{-}:= \msbar_{{g^{-}_1,m_1}}\times\ldots\times \msbar_{{g^{-}_r,m_r}}$ denotes the canonical quotient morphism. 
Note that by construction all fibres of $f^-_{S'}: {\cal C}^{-}_{S'} \rightarrow S'$ are isomorphic to $C_0^{-}$, hence ${\theta}'$ is constant. So there is in fact a commutative diagram
\[ \begin{array}{rcccc}
E^{-}_{S'} &\stackrel{\theta_{S'}^{-}}{\longrightarrow} &F&\subset& \hilb{g^{-},n,m}^{-} \\[3pt]
\mbox{\tiny ${p_{S'}^{-}}$} \downarrow \:\:&&\downarrow  &&\downarrow \mbox{\footnotesize ${\pi}$}\\[1pt]
S'\: &\longrightarrow & \{ [C_0^{-}]\}&\subset&\msbar_{{g^{-},m}}^{-}
\end{array}\]
where $F\cong P^{-} / A$ is the fibre of the canonical
morphism $\hilb{g^{-},n,m}^{-} \rightarrow \msbar_{{g^{-},m}}^{-}$
over the point $[C_0^{-}]$ representing the isomorphism class of $C_0^-$.

Fix an embedding of $C_0^{-}$ into $\PP^{N_1^{-}}\times \ldots\times
\PP^{N_r^{-}}$, corresponding to a point in $\hilb{g^{-},n,m}^{-}$,
which shall also be denoted by $[C_0^{-}]$. Note that the decomposition of $A =
\Aut(C_0^{-})$ into a product $\Aut(C_1^-)\times\ldots\times \Aut(C_r^-)$ is compatible with the decomposition of $P^{-}$, since each connected component of $C_0^-$ contains at least one marked point.  
A standard computation now shows that 
\[ E^{-}_{S'}/A \cong F \times S' \]
is a trivial fibration. Indeed, we can define a trivializing section
$\sigma': S' \rightarrow E^{-}_{S'}/ A$  as follows.
For $s \in S'$ choose some $e\in (p_{S'}^{-})^{-1}(s) \subset
E_{S'}^{-}$, such that $\theta_{S'}^{-}(e) = [C_0^{-}]\in F
\subset \hilb{g^{-},n,m}^{-} $. Under the identification $F \cong
P^{-}/A$ the point  $[C_0^{-}]$ corresponds to the class $[\id]$ of
the identiy $\id \in P^{-}$. Hence by the $P^{-}$-equivariance of $
\theta_{S'}^{-}$, the point $e$ is uniquely defined up to the action of
$A$, and thus $\sigma'(s) := [e]\in E_{S'}^{-}/A$ is well defined. 

$(iii)$ So far we have been working on an \'etale cover $u:S'\rightarrow S$. 
From the construction of the section $\sigma': S' \rightarrow E^{-}_{S'}/ A$ one sees that, if $s_1,s_2 \in S'$ are closed points such that $u(s_1)=u(s_2)$, then $\sigma'(s_1)$ and $\sigma'(s_2)$ differ at most by an element in $\Aut(C_0^-;\CQ)/A = \Gamma$. Thus we can define a section $ \overline{\sigma}:  S \rightarrow  E_{S'}^{-}/\Aut(C_0^-;\CQ)$. Using the inclusion $E_{S'}^-/\Aut(C_0;\calq)  \hookrightarrow  E_{S'}  / A$, and composing it with the canonical quotient morphism, as well as the pullback morphism $E_{S'}\rightarrow E$, we finally obtain a section ${\sigma} : \:\: S \rightarrow E / P\times  A $. We now define a subbundle $E^+$ of the principal
$\PGL(N+1)$-bundle $p : E \rightarrow S$ by the following Cartesian diagram:
\[ \diagram
E^+ \rto \dto_{p^+} & E \dto\\ 
S \rto_\sigma & E/P\times A .
\enddiagram\]
It is easy to see that  $ p^{+}: E^{+} \rightarrow S$ is a principal $P\times A$-bundle, and its extension
$(E^{+}\times \PGL(N+1))/( P\times  A)$ is isomorphic to $E$. 
Note that this property makes this construction the inverse of the
construction of proposition \ref{MainProp}. 

$(iv)$ 
Recall that our aim is to construct an object
$(E^{+},p^{+},\phi^{+})\in $ \linebreak $[\hgammax/P\times
A](S)$. To do this, we stil need to define an $P\times
A$-equivariant morphism $\phi^{+}:E^{+} \rightarrow
\hgammax$. 

Let $v: K' \rightarrow \kpx$ be an \'etale cover as in lemma \ref{315}, such that the pullback  $f_{K'}: {\cal C}_{K'}\rightarrow K'$ of the universal curve ${\cal C}_{g,n,0}\rightarrow \kpx $ contains
a trivial subcurve $f_{K'}^{-} : C_0^{-}\times K' \rightarrow
K'$. Let $f_{K'}^{+} : {\cal C}_{K'}^{+} \rightarrow K'$
denote the closure of the complement of this trivial subcurve.  
This is an $m$-pointed stable curve of genus $g^{+}$ over $K'$. 
By construction, $E^{+}
\subset E$ is a subbundle, which is stabilized under the action of the subgroup $P$ in $\PGL(N+1)$. Thus the $\PGL(N+1)$-equivariant morphism $\phi: E\rightarrow \kpx$ restricts to a
$P$-equivariant morphism
${\phi}_{K}^+ :    E^{+} \rightarrow \kpx$.

Consider the pullback $ \phi_{v}^+:  E^+_v := v^\ast E^+  \rightarrow K'$ of this morphisms along $v$.  We define an $m$-pointed stable curve $f_v^{+}: {\cal C}_v^{+} \rightarrow E^{+}_v$ of genus $g^+$ as the pullback of the stable curve $f_{K'}^{+} : {\cal C}_{K'}^{+} \rightarrow K'$ anlong $ \phi_{v}^+$. 
The projective bundle 
$ \PP(f_{v}^{+})_\ast(\omega_{{\cal C}_{v}^{+}/E^{+}_v}(S_1+\ldots+S_m))^{\otimes n}  $ 
on $E_v^{+}$ is isomorphic to the pullback of the projective bundle
$ \PP(f_{K'}^{+})_\ast (\omega_{{\cal C}^{+}_{K'}/K'}({\cal S}_1+\ldots+{\cal S}_m))^{\otimes n} $
on $K'$ along the morphism ${\phi}_v^+$. Here again $S_1,\ldots,S_m$ and ${\cal S}_1,\ldots,{\cal S}_m$ denote the divisors of the $m$ marked points on ${\cal C}_{v}^{+}$ and ${\cal C}_{K'}^{+}$, respectively. 

$(v)$ 
Note that there is a distinguished trivialization of the projective
bundle  $\PP(f_{K'}^{+})_\ast (\omega_{{\cal
C}^{+}_{K'}/K'}({\cal S}_1+\ldots+{\cal S}_m))^{\otimes n}$
over $K'$.\footnote{In fact, the triviality of the pullback $\tilde{\phi}_\beta^\ast \PP(u_\beta^{+})_\ast (\omega_{{\cal
C}^{+}_\beta/K_\beta}({\cal S}_1+\ldots+{\cal S}_m))^{\otimes n}$ shows that the principal $P$-bundle $E_\beta^+/A$ is isomorphic to the principal $P$-bundle associated to the projective bundle $\PP(u_\beta^{+})_\ast (\omega_{{\cal
C}^{+}_\beta/K_\beta}({\cal S}_1+\ldots+{\cal S}_m))^{\otimes n}$ over $K_\beta$.} 
This can be seen as follows. The curve $f_{K'}: {\cal C}_{K'}\rightarrow K'$ has been constructed as the pullback of the universal curve on $\overline{H}_{g,n,0}$, restricted to $\kpx$. Therefore there exists a natural embedding of the subcurve $f_{K'}^+: {\cal C}_{K'}^+\rightarrow K'$ into $\PP^N \times K'\rightarrow K'$. To distinguish a trivialization of $\PP(f_{K'}^{+})_\ast (\omega_{{\cal
C}^{+}_{K'}/K'}({\cal S}_1+\ldots+{\cal S}_m))^{\otimes n}$ over $K'$, it suffices to distinguish a corresponding embedding of the subcurve $f_{K'}^+: {\cal C}_{K'}^+\rightarrow K'$ into $ \PP^{N^+} \times K'\rightarrow K'$. 

Consider first the curve $f^\times: {\cal C}_{g^+,n,m}^\times \rightarrow \hx$, which is the restriction of the universal curve on the full Hilbert scheme. By definition, this curve is embedded into $ \PP^{N^+}\times \hx = \PP f^\times_\ast(\omega_{{\cal C}_{g^+,n,m}^\times/\hx}({\cal S}_1+\ldots+{\cal S}_m))^{\otimes n}$, which in turn is embedded into $\PP(u_0)_\ast(\omega_{{\cal C}_0/\hx})^{\otimes n} = \PP^N \times \hx$ by the natural inclusion of
$ \PP f^\times_\ast(\omega_{{\cal C}_{g^+,n,m}^\times/\hx}({\cal S}_1+\ldots+{\cal S}_m))^{\otimes n} $ in $
\PP(u_0)_\ast(\omega_{{\cal C}_0/\hx})^{\otimes n} $. 
Compare construction \ref{const_part_1} for details and for the notation. Note that the embedding is fibrewise independent of the ordering of the labels of the marked points. We used this in proposition \ref{const_prop} for the construction of the  morphism $\thx: \hx \rightarrow \overline{H}_{g,n,0}$. Recall that $\thx$ induces a morphism $\Theta: \hgammax \rightarrow \overline{H}_{g,n,0}$, see proposition \ref{312}. 

For each closed point $[C]\in \Theta(\hgammax) = \thx(\hx)$ representing a stable curve $C$ embedded into $\PP^N$, there is a unique subcurve isomorphic to $C_0^-$, and thus a unique subcurve $C^+$, which is an $\mgamma$-pointed stable curve of genus $g^+$, satisfying $\Theta([C^+]) = [C]$. Thus, the restricition of $f^+_{K'}:  {\cal C}_{K'}^+\rightarrow K'$ to the preimage $T'$ of $\Theta(\hgammax)$ has a distinguished embedding into $ \PP^{N^+}\times T' \rightarrow T'$, namely the one given by the embedding of the universal curve over $\hx$. We want to extend this embedding to all of $K'$, using lemma \ref{314}, which implies that $\kpx = \Theta(\hgammax) \cdot \PGL(N+1)$. Note that if $x,y\in K'$ are closed points with $v(x)=v(y)$, then the fibres of $f_{K'}^+$ over the two points differ only by a reordering of the labels of their marked points, which does not affect their embedding into $\PP^{N^+}$. 

Let $[C]\in T$ be a closed point representing a stable curve $C$ embedded into $\PP^N$, with a decomposition $C= C^+ \cup C_0^-$ as above, together with the embedding of $C^+$ into $\PP^{N^+}$. For an element $\gamma\in \PGL(N+1)$, the distinguished $\mgamma$-pointed stable subcurve of genus $g^+$ of the translated embedded curve $\gamma C$ is the curve $\gamma C^+$, by the uniqueness statement of corollary \ref{K324}. The obvious embedding of $\gamma C^+$ into $\gamma \PP^{N^+} \subset \PP^N$ is in fact independent of the choice of $\gamma$. 

Indeed, let $\gamma'\in \PGL(N+1)$ be another element such that $\gamma'([C])=\gamma([C])$. Then $\gamma'\circ\gamma^{-1}\in\Aut(C)$. By the definition of $\cdx$, this group splits as
$\Aut(C)=\Aut_\Gamma(C^+;\calq)\times A$. Since the embedding of $C^+$ into $\PP^{N^+}$ is independent of the ordering of the labels of the marked points, and since $\Aut(C^+) = \Stab_{\PGL(N^+ +1)}([C^+])$, any elements of $\Aut(C)$ lying in $\Aut_\Gamma(C^+;\calq)$ leave the embedding invariant. By construction \ref{const_part_2}, there is a fixed embedding $C_0^- \subset \PP^{N^-} \subset \PP^N$, and thus  inclusions $A\subset \PGL(N^- +1) \subset \PGL(N+1)$. Note that $\PGL(N^- +1)$ acts trivially on $\PP^{N^+}$ outside the intersection $\PP^{N^+}\cap \PP^{N^-}$. However, elements of $A$ fix by definition the marked points of $C_0^-$, which span the intersection $\PP^{N^+}\cap \PP^{N^-}$ as a subspace of $\PP^N$, and thus fix the subspace $\PP^{N^+}$. Hence elements of $\Aut(C)$ lying in $A$ fix the embedding of $C^+$ as well. 

Therefore, for all closed points $[C]\in \kpx$, representing an embedded stable curve $C\subset \PP^N$, there is a well-defined embedding of the distinguished $\mgamma$-pointed stable subcurve $C^+ \subset C$ into $\PP^{N^+}$, and so there exists an embedding of  $f_{K'}^+: {\cal C}_{K'}^+\rightarrow K'$ into $\PP^{N^+} \times K'\rightarrow K'$, as claimed.  

$(vi)$ 
So finally we obtain an $m$-pointed stable curve $f_{v}^{+}: {\cal C}_{v}^{+} \rightarrow E^{+}_v$, together with a trivialization of the projective bundle $\PP(f_{v}^{+})_\ast(\omega_{{\cal C}_{v}^{+}/E^{+}_v}(S_1+\ldots+S_m))^{\otimes n}$. By the universal property of the Hilbert scheme, this defines a $P$-equivariant morphism
$ \tilde{\phi}^{+}_v : \: \: E_v^{+} \rightarrow \hx$. 

Recall that $E_v^+ = v^\ast E^+$. If $e_1,e_2\in E$ get mapped to the same point in $E^+$, then by the construction of $ \tilde{\phi}^{+}_v$ there exists a permutation $\gamma  \in\Gamma$, such that $\gamma \tilde{\phi}^{+}_v(e_1) = \tilde{\phi}^{+}_v(e_2)$. Therefore there is a well-defined induced morphism
\[ {\phi}^{+} : \quad {E}^{+} \rightarrow \hgammax,\]
which is $P\times  A$-equivariant, with $ A$
acting trivially on $\hgammax$.

$(vii)$
Taking everything together, we obtain an object 
$(E^{+}, p^{+},\phi^{+}) $ of the fibre category $ [  \hgammax/P\times A ](S)$.
Note that up to the action of $A$ this triple just represents the $\mgamma$-pointed stable curve over $S$, which is the closure of the complement of the \'etale trivial curve $C_0^{-}\times S \rightarrow S$ in $f_S:{\cal C}_S\rightarrow S$, as in lemma \ref{315}. 

Up to isomorphism, this construction of the reduction $E^{+}$ of $E$
is inverse to the extension used in the proof of proposition
\ref{MainProp}. In fact one has the equality $\Xi(\Lambda(E^{+},p^+,\phi^+)) =
(E^{+},p^+,\phi^+)$, and an isomorphism $\Lambda(\Xi(E,p,\phi)) \cong (E,p,\phi)$. 

$(viii)$ It finally remains to define the inverse functor $\Xi$ for all schemes $S$ on morphisms
$ \varphi :  (E_2,p_2,\phi_2) \rightarrow (E_1,p_1,\phi_1) $
in $\cdx(S)$. Such a morphism is given by a morphism of principal $\PGL(N+1)$-bundles $f : E_1 \rightarrow E_2$, such that the diagram
\[ \diagram
E_1 \xto[rrd]^{\phi_1} \drto^f \xto[rdd]_{p_1}\\
& E_2 \rto_{\phi_2} \dto^{p_2}& \kpx\\
&S
\enddiagram \]
commutes. We need to show that $f$ restricts to a morphism of principal $P\times A$-bundles
$ f^{+} : \: \: E_1^{+} \rightarrow E_2^{+} $, 
where $E_1^{+}$ and $E_2^{+}$ are the $P\times  A$-subbundles of $E_1$ and $E_2$ constructed as above. Note that this restriction, if it exists, necessarily satisfies 
$ \phi_1 = \phi_2\circ f^+ $. 

For $i=1,2$, the morphism $\phi_i$ factors as $\phi_i=\Theta\circ\phi_i^+$, where $\Xi(E_i,p_i,\phi_i)=(E_i^+,p_i^+,\phi_i^+)$. The injectivity of $\Theta : \hgammax \rightarrow \kpx$, which follows from proposition \ref{312},  implies the identity 
$ \phi_1^+ = \phi_2^+\circ f^+$, 
so $f^+:E_1^+\rightarrow E_2^+$ defines indeed a morphism between the triples $(E_2^+,p_2^+,\phi_2^+)$ and $(E_1^+,p_1^+,\phi_1^+)$. 
It is clear from the construction that this inverts the definition of $\Lambda$ on morphisms. Thus $\Xi$ is shown to be an inverse functor to $\Lambda$, up to isomorphism, provided that we can prove the existence of $f^+$. 

Clearly, such a restriction $f^+:E_1^+\rightarrow E_2^+$ exists, if and only if the condition
$ f(E_1^{+}) \subseteq E_2^{+}$  is satisfied.

To do this, let us briefly recall our construction of $\Xi$ on objects. Let $i=1,2$.  We used the morphism $\phi_i : E_i \rightarrow \kpx$ to pull back the universal curve ${\cal C}_{g,n,0}\rightarrow \kpx$ to $E_i$. This pullback contains a distinguished $m$-pointed prestable curve $f_{E_i}^- : {\cal C}_{E_i}^- \rightarrow E_i$ as a subcurve, which in turn induces a morphism $\theta_i^{-} :   E_i^{-} \rightarrow \hilb{g^{-},n,m}^{-}$. 
Strictly speaking, this is true only up to \'etale coverings. However, the above inclusion condition is not affected by \'etale pullbacks, so we may safely ignore this technical complication.

By the universal property of the Hilbert space, we have the identity
\[ (\ast) \qquad \qquad \theta_1^- = \theta_2^-\circ f.    \]
We saw that $\theta_i^-$ is in fact a morphisms $ \theta_i^- : E_i^{-} \rightarrow F$,  where $F\cong P^{-}/A$ is the fibre of the quotient morphism $\hilb{g^{-},n,m}^{-} \rightarrow \msbar_{g^{-},m}^{-}$ over the point representing the isomorphism class of $C_0^{-}$. 
Again, we identify here the image  $[\id]$ of the unit $\id \in
P^{-}/A$ with the point $[C_0^{-}]$ representing the fixed embedding of $C_0^{-}$. 

There is a section $\sigma_i': S \rightarrow E_i^{-}/A$, sending a point $s\in S$ to the equivalence class of points $e\in E_i^{-}$ in the fibre over $s$  with $\theta_i^-(e) = [\id]\in F$. This is well-defined because $\theta_i^-$ is $P^{-}$-equivariant. 
Let ${\sigma}_i$ denote the morphism from $S$ to $E_i/P\times A$ induced by $\sigma_i'$. 
Now the bundle $E_i^{+}$ is defined as the unique subbundle of $E_i$,
which fits into the Cartesian diagram
\[ \diagram
E_i^{+}\rrto \dto_{p_i^{+}} &&E_i\dto^{q_i}\\
S\rrto_{{\sigma}_i}&&E_i/P\times  A ,
\enddiagram\]
where $q_i : E_i \rightarrow E_i/P\times  A $ denotes the natural quotient map. In other words,
for an element $e\in E_i$ holds $e\in E_i^{+}$ if and only if
$ q_i (e) = {\sigma}_i(p_i^{+}(e))$, 
which in turn is equivalent to
$ \theta_i^-(e) = [\id]$
by the definition of ${\sigma}_i$. 
So the relation $f(E_1^{+}) \subseteq E_2^{+}$ translates into the condition
$\theta_2^-(f(e)) = [\id] $
for all $e\in E_1^{+}$. However, above we found the identity $(\ast) \:\: \theta_1^- = \theta_2^-\circ f$, which concludes the proof.
\ebew

\begin{cor}
An isomorphism of stacks
$ \cdx \: \cong \: \cmgammadetailx $ exists 
if and only if $\Aut(C_0^{-})$ is trivial.
\end{cor}

\proof
This follows immediately from theorem \ref{MainThm}.
\ebew

\begin{cor}\em
There exists a natural morphism of moduli stacks  
\[ \Omega : \quad \cmgammadetailx \:   \rightarrow \: \cdx, \]
which is representable, surjective, finite and unramified,  and of degree equal to the order of the group $\Aut(C_0^-)$.
\end{cor}

\proof
Using the notation from the proof of theorem \ref{MainThm}, the morphism $\Omega$ is given as the composition of the canonical quotient morphism  $[ \hgammax/P] \rightarrow [ \hgammax/P\times A]$ with the isomorphism $\Lambda$ from proposition \ref{MainProp}. Note that the quotient morphism is a representable, finite and unramified surjective morphism. Since both $\cmgammadetailx$ and $\cdx$ are smooth stacks of dimension one by lemma \ref{140}, this implies that $\Omega$ is \'etale as well. 
\ebew 

\begin{rk}\em\label{summingup}
Summarizing the above results, we obtain 
the following  commutative diagram:
\[ \begin{array}{ccccc}
\hgammax& =   & \hgammax & \hookrightarrow &\hilb{g,n,0}\\
 \downarrow &&\downarrow &&\downarrow \\
\cmgammax &\leftarrow & \cdx &\hookrightarrow & \cmgbar\\
\downarrow &&\downarrow &&\downarrow \\
\mgammax & \cong &\dx &\hookrightarrow &\msbar_{g}.
\end{array} \]
Note that, since $\Aut(C_0^-)$ acts trivially, there is also a canonical morphism
\[ \cdx \cong [ \, \hgammax / P\times A \,] \:  \longrightarrow \: [\, \hgammax / P \,] \cong \cmgammax, \]
which is of degree equal to $\frac{1}{\# \Aut(C_0^{-})}$. In \cite{Z2} it is shown that the moduli stack of $\mgamma$-pointed stable curves is a quotient of the moduli stack of $m$-pointed stable curves with respect to the action of $\Gamma$. This implies in particular that there is a finite surjective morphisms 
\[ \overline{\Omega} : \quad \cmgx \:   \rightarrow \: \cdx, \]
which is of degree equal to the order of $\Aut(C_0^-;\calq)$. In general, this morphism is not representable, see \cite{Z1}. 
\end{rk}

\section{The case of genus three}\label{Sect81}

To illustrate the  description of one-dimensional boundary substacks of the moduli space of Deligne-Mumford stable curves as given above, we want to consider the case of genus $g=3$ in detail. 

\begin{rk}\em\label{R8301}
$(i)$ There are exactly five different stable curves of genus $g=3$ with $3g-3=6$ nodes. We denote them by $C_1,\ldots,C_5$, in the same order as they are listed in Faber's paper \cite{Fa}. Schematically they can be represented as in the following pictures.
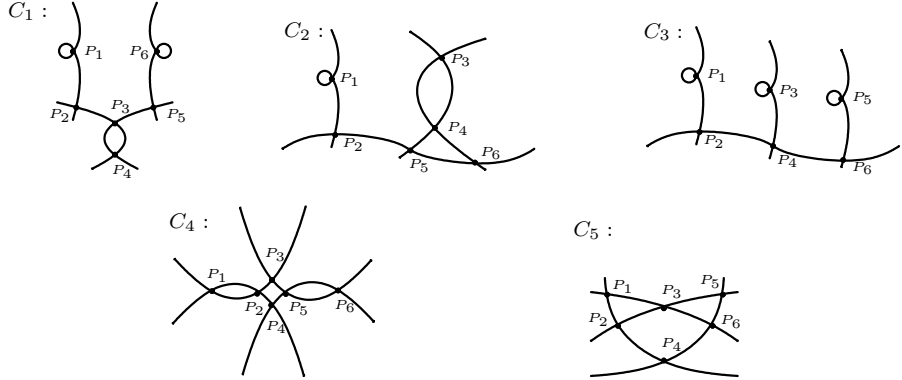
\begin{center}
\setlength{\unitlength}{0.00012in}
\begingroup\makeatletter\ifx\SetFigFont\undefined%
\gdef\SetFigFont#1#2#3#4#5{%
  \reset@font\fontsize{#1}{#2pt}%
  \fontfamily{#3}\fontseries{#4}\fontshape{#5}%
  \selectfont}%
\fi\endgroup%
{\renewcommand{\dashlinestretch}{30}

}
\end{center}
$(ii)$ 
There are eight different types of stable curves of genus $3$ with exactly $5$ nodes. Following Faber's notation, we denote them by $(a),\ldots,(h)$. They are represented by the pictures below.

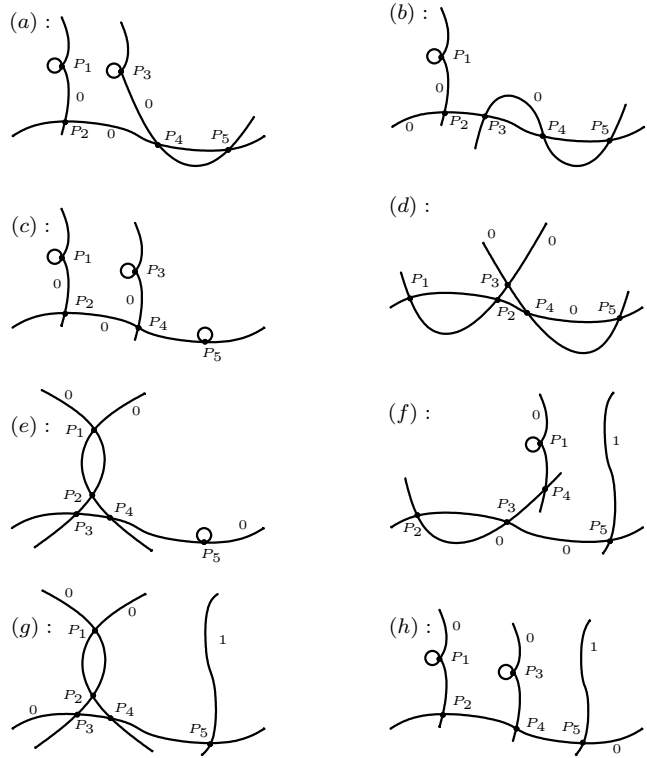
\begin{center}
\setlength{\unitlength}{0.00012in}
\begingroup\makeatletter\ifx\SetFigFont\undefined%
\gdef\SetFigFont#1#2#3#4#5{%
  \reset@font\fontsize{#1}{#2pt}%
  \fontfamily{#3}\fontseries{#4}\fontshape{#5}%
  \selectfont}%
\fi\endgroup%
{\renewcommand{\dashlinestretch}{30}
%
}
\end{center}

The number next to an irreducible component gives its geometric genus.
\end{rk}

\begin{rk}\em
The type of a curve containing $3g-4=5$ nodes does not vary within a connected component $D$ of the boundary stratum $M_3^{(5)}$. Recall that each such component can be written as ${D}= {D}(C_i \, ;P_j)$ for some $1\le i\le 5$ and  $1\le j\le 6$.  
The following table lists all possible presentations for all connected components.
\[\begin{array}{llcccc}
{D} \qquad & {D}(C_i \, ;P_j)& M_{g^+,m}'& \deg(\overline{\Omega})& \Aut(C_i^-) & \Gamma(C_i^-;\CQ)\\
\hline\\ 
(a)&D(C_1;P_{2,5})& M_{1,2}'& 4&\ZZ_2&\Sigma_2\\
   &D(C_2;P_{3,4})& M_{1,2}'& 4&\ZZ_2&\Sigma_2\\
(b)&D(C_2;P_{5,6})& M_{0,4}&12&\ZZ_2& \Sigma_3\\
(c)&D(C_1;P_{3,4})& M_{1,2}'&8&\ZZ_2\times\ZZ_2&\Sigma_2\\
   &D(C_3;P_{2,4,6})& M_{1,2}'&8&\ZZ_2\times\ZZ_2&\Sigma_2\\
(d)&D(C_4;P_{3,4})& M_{0,4}&8&\id&D_8\\
   &D(C_5;P_{1,\,\ldots,6})& M_{0,4}&8&\id&D_8\\
(e)&D(C_2;P_2)    & M_{1,2}'&4&\ZZ_2&\Sigma_2\\
   &D(C_4;P_{1,2,5,6})    & M_{1,2}'&4&\ZZ_2&\Sigma_2\\
(f)&D(C_1:P_{1,6})& M_{1,1}&4&\ZZ_2\times\ZZ_2&\id\\
(g)&D(C_2;P_1)&M_{1,1}&4&\ZZ_2\times\ZZ_2&\id\\
(h)&D(C_3;P_{1,3,5})&M_{1,1}&8&D_8&\id
\end{array}\]
The first column specifies the connected component ${D}$  by indicating the type of a general point. The second column lists all possible choices of curves  $C_i$ and nodes $P_j$ to represent the respective component. More than one subscript, like for example $P_{2,5}$, indicates that both nodes $P_2$ and $P_5$ produce the same component $D$. The third column gives the moduli scheme $M_{g^+,m}'$ covering ${D}$ as in lemma \ref{decompose}, and the fourth column states the degree of the finite morphism $\overline{\Omega}:\cmgx \rightarrow {\cal D}^\times(C_i;P_j)$ as in remark \ref{summingup}. Columns five and six give the corresponding group of automorphisms $ \Aut(C_i^-)$ and the subgroup of permutations $ \Gamma(C_i^-;\CQ) \subset \Sigma_m$, respectively. Here, $D_8$ denotes the dihedral group of order $8$. 
\end{rk}

\begin{ex}\rm
Consider the stable curve $C_2$. There are two pairs of equivalent nodes on $C_2$, in the sense that for each pair there exists an automorphism of $C_2$ interchanging the two nodes. Therefore deformations, which preserve all but one node of $C_2$, distinguish four different connected components $D$ of $M_3^{(5)}$. Let us  have a closer look at all four of those components. In particular, we want to apply the isomorphism from theorem \ref{MainThm} to the corresponding moduli substacks ${\cal D}$. We quote from \cite{Z1} to supply additional information about possible extensions of this isomorphism.

$(i)$ We consider first the connected  component ${D}={D}(C_2;P_1)$. A general point of this component represents a curve of type $(g)$. Since $m=1$, all automorphism groups split. In particular, we have the equality $D^\times(C_2;P_1) = D(C_2;P_1)$. Thus  we have an isomorphism 
$ {\cal D}(C_2;P_1) \cong {\cal M}_{1,1} / {\mathbb Z}_2 \times {\mathbb Z}_2$ by theorem \ref{MainThm}. 
In fact, in this case the isomorphism extends to the respective closures $\overline{\cal D}(C_2;P_1)  \cong  \overline{\cal M}_{1,1} / {\mathbb Z}_2 \times {\mathbb Z}_2$.

$(ii)$ Consider now $D=D(C_2;P_2)$, where a general point of this component represents a curve of type $(e)$. Note that for all points $[C^+]\in M_{1,2}'$, representing a stable curve $C^+$ with two marked points, the automorphism group $\Aut(C)$  splits. Thus there is an  isomorphism    
${\cal D} (C_2;P_2)  \cong {\cal M}_{1,2/\Sigma_2}' / \ZZ_2$. 
It can be extended to a partial compactification of  ${\cal D}^\times (C_2;P_2)$ over the point $[C_2]$, but not to its  closure.

$(iii)$ The case  of  ${D}={D}(C_2;P_3)$ is analogous to case $(ii)$. Note that this is the same connected component as ${D}(C_2;P_4)$, where a general point represents a curve of type $(a)$. There is also an  isomorphism    
${\cal D} (C_2;P_3)  \cong {\cal M}_{1,2/\Sigma_2}' / \ZZ_2$, which   can be extended to a partial compactification, but not to all of the  closure.

$(iv)$ Finally, consider the connected component ${D}={D}(C_2;P_5)={D}(C_2;P_6)$, with general point of type $(b)$. Again we have an equality $D^\times(C_2;P_5) = D(C_2;P_5)$. There is an isomorphism of stacks $ {\cal D}(C_2;P_5)   \cong {{\cal M}_{0,4/\Sigma_3}} / {\mathbb Z}_2 $. 
This isomorphism cannot be extended any further. However, one can find a natural reduced substack $\hat{\cal M}_{2,1}$ of the moduli stack $\overline{\cal M}_{2,1}$ of $1$-pointed stable curves of genus $2$, so that for the closure of $ {\cal D}(C_2;P_5)$ there is an isomorphism $\overline{\cal D}(C_2;P_5)   \cong  \hat{\cal M}_{2,1} / {\mathbb Z}_2 $. For full details, see \cite{Z1}.
\end{ex}

\begin{ex}\em
To provide an example, where the isomorphism of theorem \ref{MainThm} does not hold for all of $D$, look at the connected component ${D}={D}(C_3; P_2)$. A general point represents a curve of type $(c)$. If $[C^+]\in M_{1,2}'$ represents a stable curve $C^+$ with two marked points, which are in a special position, then the exact sequence of automorphism groups of lemma \ref{newsplit} does not split.  There is an isomorphism
$ {\cal D}^\times(C_3;P_2)   \cong  {{\cal M}^\times_{1,2/\Sigma_2}} / {\mathbb Z}_2 \times {\mathbb Z}_2 $, 
but it cannot be extended any further. 
\end{ex}

\end{document}